\documentclass[11pt,a4paper]{article}
\RequirePackage{amsmath}%
\RequirePackage{amsthm}%
\usepackage{amsfonts}%
\usepackage{amssymb}%
\usepackage{graphicx}%
\usepackage{cite}
\usepackage{url}

\renewcommand{\Re}{\mathop{\mathrm{Re}}}
\renewcommand{\Im}{\mathop{\mathrm{Im}}}

\renewcommand{\i}{\mathrm{i}}

\newcommand{\diag}{\mathrm{diag}}

\renewcommand{\Re}{\mathop{\mathrm{Re}}}
\renewcommand{\Im}{\mathop{\mathrm{Im}}}
\renewcommand{\i}{\mathrm{i}}

\newcommand{\CC}{{\mathbb C}}

\newcommand{\bM}{{\bf M}}
\newcommand{\bN}{{\bf N}}
\newcommand{\bI}{{\bf I}}
\newcommand{\bJ}{{\bf J}}

\newcommand{\bs}{{\bf s}}
\newcommand{\bt}{{\bf t}}
\newcommand{\bx}{{\bf x}}
\newcommand{\bz}{{\bf z}}
\newcommand{\by}{{\bf y}}

\newcommand{\ba}{{\bf a}}
\newcommand{\bone}{{\bf 1}}

\newcommand{\bvb}{{\bf b}}

\newcommand{\bvd}{{\bf d}}
\newcommand{\bve}{{\bf e}}

\begin{document}

\title{Fast solution of boundary integral equations with the generalized Neumann kernel}

\author{Mohamed M.S. Nasser}

\date{}
\maketitle

\vskip-0.8cm %
\centerline{Department of Mathematics, Faculty of Science, King Khalid University,} %
\centerline{P. O. Box 9004, Abha 61413, Saudi Arabia.}%
\centerline{E-mail: mms\_nasser@hotmail.com}

\begin{center}
\begin{quotation}
{\noindent {\bf Abstracts.\;\;}%
A fast method for solving boundary integral equations with the generalized Neumann kernel and the adjoint generalized Neumann kernel is presented. The complexity of the presented method is $O((m+1)n\ln n)$ for the integral equation with the generalized Neumann kernel and $O((m+1)n)$ for the integral equation with the adjoint generalized Neumann kernel where $m+1$ is the multiplicity of the multiply connected domain and $n$ is the number of nodes in the discretization of each boundary component. The presented numerical results illustrate that the presented method gives accurate results even for domains with high connectivity, domains with piecewise smooth boundaries, and domains with close boundaries.
}%
\end{quotation}
\end{center}

\begin{center}
\begin{quotation}
{\noindent {\bf Keywords.\;\;}%
Generalized Neumann kernel; boundary integral equations; Nystr\"om method; Fast Multipole Method: GMRES; numerical conformal mapping.
}%
\end{quotation}
\end{center}

\begin{center}
\begin{quotation}
{\noindent {\bf MSC.\;\;} 45B05; 65R20; 30C30.}
\end{quotation}
\end{center}
 
\section{Introduction}
\label{sc:int}

This paper is concerned with:
\begin{enumerate}
	\item Numerical solution of the integral equation with the generalized Neumann kernel
\begin{equation}\label{e:ie}
(\bI-\bN)\mu=-\bM\gamma
\end{equation}
and numerical computing of the piecewise constant function $h$ given by
\begin{equation}\label{e:h}
h=[\bM\mu-(\bI-\bN)\gamma]/2.
\end{equation}
	\item Numerical solution of the integral equation with the adjoint generalized Neumann kernel
\begin{equation}\label{e:ie*}
(\bI+\bN^\ast+\bJ)\mu=\gamma.
\end{equation}
\end{enumerate}
See \S\ref{sc:op} below for the definitions of the operators $\bI$, $\bN$, $\bM$, $\bN^\ast$, and $\bJ$. Both boundary integral equations have been used to solve several problems in mathematics and mathematical physics in multiply connected domains such as the numerical conformal mapping~\cite{Nas-cmft09,Nas-siam09,Nas-apam,Nas-jmaa11,Nas-jmaa13,Nas-siam13,Nas-inv,Yun-sp,Yun-inv}, the Riemann-Hilbert problem~\cite{Nas-cr,Weg-Mur-Nas,Weg-Nas}, the Dirichlet problem~\cite{Nas-amc11,Weg-Mur-Nas}, the Neumann problem~\cite{Nas-amc11}, the mixed boundary value problem~\cite{AlH-bvp,AlH-aip,Nas-jam}, and the potential flow problem~\cite{Nas-cmft11,Nas-Sak}.

For bounded or unbounded multiply connected domains of connectivity $m+1$, discretizing the boundary integral equations~(\ref{e:ie}) and~(\ref{e:ie*}) by the Nystr\"om method with the trapezoidal rule yields dense and nonsymmetric $(m+1)n\times(m+1)n$ linear systems where $n$ is the number of nodes in the discretization of each boundary component. The rate of the convergence of the Nystr\"om method with the trapezoidal rule depends on the smoothness of the integrands which in turn depends on the smoothness of the boundary. If the boundaries are of class $C^{q+2}$ and the function $\gamma$ is of the class $C^q$, then the rate of the convergence of the Nystr\"om method with the trapezoidal rule is $O(1/n^q)$. For analytic boundaries and analytic $\gamma$, the Nystr\"om method with the trapezoidal rule converges exponentially~\cite{Kre90}. For domains with corners, accurate results are obtained if we use the trapezoidal rule with a graded mesh~\cite{Kre99,Kre90}. A fast method for solving the $(m+1)n\times(m+1)n$ linear systems obtained by discretizing the integral equation~(\ref{e:ie}) has been presented in~\cite{Nas-siam13} where the linear system is solved by the generalized minimal residual (GMRES) method. Each iteration of the GMRES method requires a matrix-vector product which can be computed using the Fast Multipole Method (FMM) in $O((m+1)n)$ operations. However, the discretization of the singular operator $\bM$ in~(\ref{e:ie}) and~(\ref{e:h}) requires  $O((m+1)n\ln n)$ operations. Discretizing the operator $\bM$ in~\cite{Nas-siam13} requires computing the derivatives $\gamma'$ and $\mu'$ where the derivative of the known function $\gamma$ can be computed analytically and the derivative of the unknown function $\mu$ should be computed numerically. 

This paper presents a new method for fast computing of the functions $\bM\gamma$ in~(\ref{e:ie}) and  $\bM\mu$ in~(\ref{e:h}). We assume that the functions $\gamma$ and $\mu$ are only H\"older continuous functions without any differentiability requirement. Thus, the stability issue of the numerical differentiation of the function $\mu$ is avoided. We shall rewrite the discretizing matrix of the operator $\bM$ as a sum of two matrices. The multiplication of the first matrix by a vector can be computed by the FMM in $O((m+1)n)$ operations. The second matrix is a block of $m+1$ circulant matrices, Hence, the multiplication of the second matrix by a vector can be computed by the FFT in $O((m+1)n\ln n)$ operations. Then, as in~\cite{Nas-siam13}, the discretized linear system is solved by a combination of the GMRES method and FMM in $O((m+1)n)$ operations. Hence, the unique solution $\mu$ of the integral equation~(\ref{e:ie}) and the $h$ in~(\ref{e:h}) are computed in $O((m+1)n\ln n)$ operations. This paper presents also a new method for fast solution of the integral equation with the adjoint generalized Neumann kernel~(\ref{e:ie*}) in $O((m+1)n)$ operations. Based on the presented methods, two MATLAB functions will be presented in this paper:
\begin{enumerate}
	\item \verb|FBIE|: for fast solving the integral equation with the generalized Neumann kernel~(\ref{e:ie}) and fast computing the piecewise constant function $h$ in~(\ref{e:h}).
	\item \verb|FBIEad|: for fast solving the integral equation with the adjoint generalized Neumann kernel~(\ref{e:ie*}).
\end{enumerate}

The solutions of the integral equations~(\ref{e:ie}) and~(\ref{e:ie*}) yield the boundary values of the conformal mapping and the solutions of the boundary value problems. Computing the interior values required computing the Cauchy integral formula. For the convenience of the reader, we present a MATLAB function \verb|FCAU| for fast computing of the Cauchy integral formula. The presented MATLAB functions, \verb|FBIE|, \verb|FBIEad|, and \verb|FCAU|, will be useful for computing the conformal mapping and solving potential flow problems of domains of high connectivity, see e.g.,~\cite{Nas-Sak}. 

The performance of the presented method has been tested on four numerical examples which include domains with high connectivity, domains with piecewise smooth boundaries, and domains with close boundaries. In the first example, we consider a bounded and unbounded domain of connectivity more than one thousands. For the bounded domain, half of the boundaries are piecewise smooth. In the second, third and fourth examples, we consider multiply connected domains with close boundaries. The distance between the boundaries can be as small as $10^{-4}$.

Other integral equations which has been solved by the FMM are the potential theory boundary integral equations~\cite{Gre93,Hel-Oja,Hel-Wad,Rok85} and the Kerzman-Stein integral equation~\cite{Don-Rok}. For more details on FMM and GMRES, see~\cite{Che05,Gre-Rok,Liu09,Rok85,Saa-Sch}.

\section{Auxiliary material}
\label{sc:aux}
\subsection{The Multiply Connected domain}

\begin{figure}%
\centerline{\scalebox{0.6}{\includegraphics{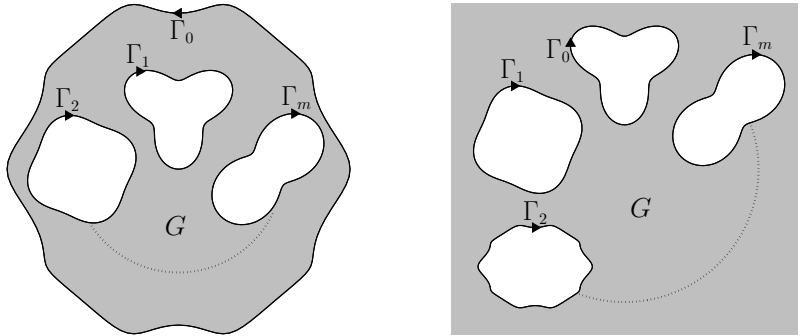}}}
\caption{\rm The bounded (left) and unbounded (right) multiply connected domain $G$ of connectivity $m+1$.} 
\label{f:dom}
\end{figure}

Let $G$ be an $(m+1)-$multiply connected domain in the extended complex plane $\overline{\CC}:=\CC\cup\{\infty\}$. The domain $G$ can be bounded or unbounded. For bounded $G$, we assume that $\alpha$ is a fixed point in $G$. If $G$ is unbounded, then $\infty\in G$. Let $G$ has the boundary
\[
\Gamma:=\partial G= \cup_{j=0}^{m} \Gamma_j
\]
where $\Gamma_0, \Gamma_1, \ldots, \Gamma_m$ are closed Jordan curves. 
The orientation of $\Gamma$ is such that $G$ 
is always on the left of~$\Gamma$. See Fig.~\ref{f:dom}. 

For $j=0,1,\ldots,m$, the curve $\Gamma_j$ is parametrized by a $2\pi$-periodic twice continuously differentiable complex function $\eta_j(t)$ with non-vanishing first derivative $\eta'_j(t)\ne 0$ for $t\in J_j:=[0,2\pi]$. The total parameter domain $J$ is the disjoint union of $m+1$ intervals $J_0,J_1,\ldots,J_m$, 
\begin{equation}\label{e:J}
J = \bigsqcup_{j=0}^{m} J_j=\bigcup_{j=0}^{m}\{(t,j)\;:\;t\in J_j\},
\end{equation}
i.e., the elements of $J$ are order pairs $(t,j)$ where $j$ is an auxiliary index indicating which of the intervals contains the point $t$~\cite[p.~394]{Lee00}. 
We define a parametrization of the whole boundary $\Gamma$ as the complex function $\eta$ defined on $J$ by
\begin{equation}\label{e:eta-1}
\eta(t,j)=\eta_j(t), \quad t\in J_j,\quad j=0,1,\ldots,m.
\end{equation}
In this paper, we shall assume for a given $t$ that the auxiliary index $j$ is known so we replace the pair $(t,j)$ in the left-hand side of~(\ref{e:eta-1}) by $t$. Thus, the function $\eta$ in~(\ref{e:eta-1}) is written as
\begin{equation}\label{e:eta}
\eta(t):= \left\{ \begin{array}{l@{\hspace{0.5cm}}l}
\eta_0(t),&t\in J_0,\\
\eta_1(t),&t\in J_1,\\
\hspace{0.3cm}\vdots\\
\eta_m(t),&t\in J_m.
\end{array}
\right.
\end{equation}

Let $H$ be the space of all real H\"older continuous $2\pi$-periodic functions $\phi(t)$ 
of the parameter $t$ on $J_j$ for $j=0,\ldots, m$, i.e.,
\[
\phi(t) = \left\{
\begin{array}{l@{\hspace{0.5cm}}l}
 \phi_0(t),     & t\in J_0, \\
 \phi_1(t),     & t\in J_1, \\
  \vdots       & \\
 \phi_m(t),     & t\in J_m, \\
\end{array}
\right.
\]
with real H\"older continuous $2\pi$-periodic functions $\phi_0,\ldots,\phi_m$. In view of the smoothness of $\eta$, a real H\"older continuous function $\hat\phi$ 
on $\Gamma$ can be interpreted via $\phi(t):=\hat\phi(\eta(t))$, $t\in J$, 
as a function $\phi\in H$; and vice versa. The subspace of $H$ that consists of real piecewise constant functions of the form
\[
h(t) = \left\{
\begin{array}{l@{\hspace{0.5cm}}l}
 h_0,     & t\in J_0, \\
 h_1,     & t\in J_1, \\
  \vdots       & \\
 h_m,     & t\in J_m, \\
\end{array}%
\right.
\]
with real constants $h_0,h_1,\ldots,h_m$ is denoted by $S$. For simplicity, the piecewise constant function $h$ will be written as
\[
h(t)=(h_0,h_1,\ldots,h_m).
\]

\subsection{The generalized Neumann kernel}

Let $\theta$ be the piecewise constant function
\begin{equation}\label{e:thet}
\theta(t)=(\theta_0,\theta_1,\theta_2,\ldots,\theta_m)
\end{equation}
where $\theta_0,\theta_1,\theta_2,\ldots,\theta_m$ are given real constants. We define a complex-valued function $A$ on $\Gamma$ by
\begin{equation}\label{e:A}
A(t) = \left\{
\begin{array}{l@{\hspace{0.5cm}}l}
 e^{\i\left(\frac{\pi}{2}-\theta(t)\right)}\,(\eta(t)-\alpha),     & \text{if $G$ is bounded}, \\[6pt]
 e^{\i\left(\frac{\pi}{2}-\theta(t)\right)},     & \text{if $G$ is unbounded}. \\
\end{array}%
\right.
\end{equation}
The adjoint of the function $A$ is defied by
\begin{equation}\label{e:A*}
\tilde A(t) = \frac{\eta'(t)}{A(t)}.
\end{equation}

The generalized Neumann kernel formed with $A$ and $\eta$ is defined by
\begin{equation}\label{e:N}
 N(s,t) :=  \frac{1}{\pi}\Im\left(
 \frac{A(s)}{A(t)}\frac{\eta'(t)}{\eta(t)-\eta(s)}\right).
\end{equation}
We define also a kernel
\begin{equation}\label{e:M}
 M(s,t) :=  \frac{1}{\pi}\Re\left(
 \frac{A(s)}{A(t)}\frac{\eta'(t)}{\eta(t)-\eta(s)}\right).
\end{equation}
The kernel $N$ is continuous with 
\begin{equation}\label{e:N-tt}
N(t,t)= \frac{1}{\pi} \left(\frac{1}{2}\Im
\frac{\eta''(t)}{\eta'(t)} -\Im\frac{A'(t)}{ A(t)}
\right).
\end{equation}
The kernel $M$ is singular. When $s,t\in J_j$  are in the same parameter interval $J_j$, then 
\begin{equation}\label{e:M-M1}
M(s,t)= -\frac{1}{2\pi} \cot \frac{s-t}{2} + M_1(s,t)
\end{equation}
with a continuous kernel $M_1$ which takes on the diagonal the values
\begin{equation}\label{e:M1-tt}
M_1(t,t)= \frac{1}{\pi} \left(\frac{1}{2}\Re
\frac{\eta''(t)}{\eta'(t)} -\Re \frac{A'(t)}{ A(t)}
\right).
\end{equation}

The generalized Neumann kernel 
formed with $\tilde A$ and $\eta$ 
is defined by
\begin{equation}\label{e:Nt}
 \tilde N(s,t) =  \frac{1}{\pi}\Im\left(
 \frac{\tilde A(s)}{\tilde A(t)}\frac{\eta'(t)}{\eta(t)-\eta(s)}\right).
\end{equation}
Using the definition of the adjoint function $\tilde A$, we have
\begin{equation}\label{e:Nt-N*}
\tilde N(s,t)=-N(t,s)=-N^\ast(s,t)
\end{equation}
where $N^\ast$ is the adjoint of the generalized Neumann 
kernel $N$.
Similarly, the kernel $\tilde M$ defined by
\begin{equation}\label{e:Mt}
 \tilde M(s,t) =  \frac{1}{\pi}\Re\left(
 \frac{\tilde A(s)}{\tilde A(t)}\frac{\eta'(t)}{\eta(t)-\eta(s)}\right)
\end{equation}
satisfies
\begin{equation}\label{e:Mt-M*}
 \tilde M(s,t)=-M(t,s) = -M^\ast(s,t).
\end{equation}

For more details on generalized Neumann kernel, see~\cite{Weg-Mur-Nas,Weg-Nas}.

\subsection{The integral operators}
\label{sc:op}

The integral operator
\begin{equation}\label{e:bN}
  \bN \mu(s) := \int_J N(s,t) \mu(t) dt, \quad s\in J,
\end{equation}
is a Fredholm integral operator. The operator
\begin{equation}\label{e:bM}
  \bM\mu(s) := \int_J  M(s,t) \mu(t) dt, \quad s\in J,
\end{equation}
is a singular integral operator. 
The adjoint operator $\bN^\ast$ is defined by
\begin{equation}\label{e:bN*}
\bN^\ast\mu(s) := \int_J N^\ast(s,t) \mu(t) dt, \quad s\in J.
\end{equation}
We define also an integral operator $\bJ$ by
\begin{equation}\label{e:bJ}
\bJ\mu(s) := \int_J \delta(s,t) \mu(t) dt, \quad s\in J,
\end{equation}
where the kernel $\delta(s,t)$ is defined for $s\in J_k$ and $t\in J_j$, $k,j=0,1,\ldots,m$, by
\begin{equation}\label{e:delta}
\delta(s,t)=\left\{
\begin{array}{ll}
\frac{1}{2\pi}, & k=j, \\	[6pt]
0,              & k\ne j. \\
\end{array}
\right.
\end{equation}
Hence,
\begin{equation}\label{e:bJ2}
  \bJ\mu(s) = \left(\frac{1}{2\pi}\int_{J_0}\mu(t)dt,\frac{1}{2\pi}\int_{J_1}\mu(t) dt,
	\ldots,\frac{1}{2\pi}\int_{J_m}\mu(t) dt\right),
\end{equation}
i.e., the function $\bJ\mu(s)$ is a piecewise constant function.

\subsection{The trapezoidal rule}

Let $n$ be a given even positive integer. For $k=0,1,\ldots,m$, we define in each
interval $J_k$ the $n$ equidistant nodes 
\[
s_{k,p}=(p-1)\frac{2\pi}{n}\in J_k, \quad p=1,2,\ldots,n.
\]
The total number of nodes $s_{k,p}$ in the total parameter domain $J$ is $(m+1)n$. 
We shall denote these $(m+1)n$ 
nodes by $t_i$, $i=1,2,\ldots,(m+1)n$, i.e.,
\begin{equation}\label{e:ti-skp}
t_{kn+p}=s_{k,p}\in J, \quad k=0,1,\ldots,m, \quad p=1,2,\ldots,n.
\end{equation}
We define the $(m+1)n\times1$ vector $\bt$ by
\[
\bt=(t_{1},t_{2},\ldots,t_{(m+1)n})^T
\]
where $T$ denotes transportation. For any function $\gamma(t)$ defined on $J$, 
we define $\gamma(\bt)$ as the $(m+1)n
\times1$ vector obtained by componentwise evaluation of the function $\gamma(t)$ 
at the points $t_i$, $i=1,2,\ldots,(m
+1)n$. 

As in MATLAB, for any two vectors $\bx$ and $\by$, we define $\bx.*\by$ as the
componentwise vector product of $\bx$ and $\by$. 
If $\by_j\ne0$ for all $j=1,2\ldots,(m+1)n$, we define $\bx./\by$ as the
componentwise vector division of $\bx$ by $\by$.
For simplicity, we denote $\bx.*\by$ by $\bx\by$ and $\bx./\by$ by 
$\frac{\bx}{\by}$.

Let $\gamma\in H$ be a given function, i.e., $\gamma(t)$ is $2\pi$-periodic in each interval $J_k$, $k=0,1,\ldots,m$. Thus, the trapezoidal rule becomes
\begin{equation}\label{e:trap}
\int_J \gamma(t)dt =\sum_{k=0}^{m} \int_{J_k}\gamma(t)dt 
\approx \frac{2\pi}{n}\sum_{k=0}^{m}\sum_{p=1}^{n}\gamma(s_{k,p}) 
=\displaystyle\frac{2\pi}{n}\sum_{j=1}^{(m+1)n}\gamma(t_j). 
\end{equation}

\subsection{The MATLAB function: {\tt zfmm2dpart}}

In this paper we shall use function {\tt zfmm2dpart} in the MATLAB toolbox FMMLIB2D~\cite{Gre-Gim12} to compute complex-valued sums of the form
\begin{equation}\label{e:fmm-1}
\sum_{\begin{subarray}{c} j=1\\j\ne i\end{subarray}}^{(m+1)n} 
\frac{1}{\eta(t_i)-\eta(t_j)}x_j, 
\quad i=1,2,\ldots,(m+1)n, 
\end{equation}
where $x_j$ are real or complex constants. Let $\bx$ be the $(m+1)n\times1$ vector 
\begin{equation}\label{e:fmm-r}
\bx=(x_1,x_2,\ldots,x_{(m+1)n})^T
\end{equation}
and $E$ be the $(m+1)n\times(m+1)n$ matrix with the elements
\begin{equation}\label{e:fmm-E}
(E)_{ij}:= \left\{
\begin{array}{l@{\hspace{1cm}}l}
\displaystyle 0,                             & i=j,     \\[0.00cm]
\displaystyle\frac{1}{\eta(t_i)-\eta(t_j)}, & i\ne j, \quad i,j=1,2,\ldots,(m+1)n.   \\
\end{array}\right.
\end{equation}
Hence Eq.~(\ref{e:fmm-1}) can be written as a matrix-vector product $E\bx$. 
Let $\ba$ be the $2\times(m+1)n$ real vector
\begin{equation}\label{e:fmm-a}
\ba= \left(
\begin{array}{l}
	\Re\eta(\bt)^T\\[0.25cm]
	\Im\eta(\bt)^T
\end{array}
\right).
\end{equation}
Hence the matrix-vector product $E\bx$ can be computed using the MATLAB 
function {\tt zfmm2dpart} in $O((m+1)n)$ operations by
\begin{equation}\label{e:E-fmm}
E\bx = {\tt zfmm2dpart  }(iprec,(m+1)n,\ba,\bx^T,1) 
\end{equation}
where the tolerance of the FMM 
is $0.5\times 10^{-3}$ for $iprec=1$, $0.5\times 10^{-6}$ for $iprec=2$, $0.5\times 10^{-9}$ for $iprec=3$, $0.5\times 10^{-12}$ for $iprec=4$ and $0.5\times 10^{-15}$ for $iprec=5$.

Similarly, the MATLAB function {\tt zfmm2dpart} can be used to compute complex-valued sums of the form 
\begin{equation}\label{e:fmm-2}
\sum_{j=1}^{(m+1)n} \frac{1}{z_i-\eta(t_j)}x_j, \quad i=1,2,\ldots,\hat n,
\end{equation}
where $x_j$ are real or complex constants and $z_i$ are $\hat n$ given points in $G$ with a given positive integer $\hat n$. Let $F$ be the $\hat n\times(m+1)n$ matrix with the elements
\begin{equation}\label{e:F}
(F)_{ij}:= \frac{1}{z_i-\eta(t_j)},
\quad i=1,2,\ldots,\hat n,
\quad j=1,2,\ldots,(m+1)n.
\end{equation}
Let also $\bx$ be the $(m+1)n\times1$ vector $\bx=(x_1,x_2,\ldots,x_{(m+1)n})^T$, $\bz$ be the $\hat n\times1$ complex vector $\bz=(z_1,z_2,\ldots,z_{\hat n})$ and $\bvd$ be the $2\times(m+1)n$ real vector
\begin{equation}\label{e:bvc}
\bvd= \left(
\begin{array}{l}
	\Re\bz\\[0.25cm]
	\Im\bz
\end{array}
\right).
\end{equation}
Hence Eq.~(\ref{e:fmm-2}) can be written as a matrix-vector product $F\bx$. The matrix-vector product $F\bx$ can be computed using the MATLAB function {\tt zfmm2dpart} in  $O((m+1)n+\hat n)$ operations by
\begin{equation}\label{e:F-fmm}
F\bx = {\tt zfmm2dpart}(iprec,(m+1)n,\ba,\bx^T,0,0,0,\hat n,\bvd,1,0,0).
\end{equation}

\section{Solving the integral equation with the generalized Neumann 
kernel}

\subsection{The integral equation}

We shall use singularity subtraction to rewrite the operators $\bN$ and $\bM$ to make these operators more suitable for using the FMM. This procedure is useful for solving the integral equation~(\ref{e:ie}) for domains with corners (see~\cite{Ans81,Kre99,Kre90,Nas-siam13,Nas-cr,Rat93}). It is also useful for solving the integral equation~(\ref{e:ie}) for domains with close boundaries (see the numerical examples below). 

It is known that the constant function is an eigenfunction of the generalized Neumann kernel $N$ corresponding to the eigenvalue $\lambda=-1$and an eigenfunction of the singular kernel $M$ corresponding to the eigenvalue $\lambda=0$~\cite{Nas-cmft09,Nas-amc11}, i.e.,
\begin{equation}\label{e:NM-const}
\int_{J}N(s,t)dt=-1, \quad \int_{J}M(s,t)dt=0.
\end{equation}
Thus, the integral equation~(\ref{e:ie}) can be written as 
\begin{equation}\label{e:ie-m}
2\mu(s)-\int_{J}N(s,t)[\mu(t)-\mu(s)]dt = -\phi(s), 
\end{equation}
where
\begin{equation}\label{e:ie-phi}
\phi(s) =\int_{J}M(s,t)[\gamma(t)-\gamma(s)]dt.
\end{equation}
The integral equation~(\ref{e:ie-m}) is valid even if the boundary $\Gamma$ 
is piecewise smooth (see~\cite{Nas-cr}).

\subsection{The Nystr\"om method}

Discretizing the integral in~(\ref{e:ie-m}) by the 
trapezoidal rule~(\ref{e:trap}) 
and substituting $s=t_i$, we obtain the linear system
\begin{equation}\label{e:ie-sys}
2\mu(t_i)-\frac{2\pi}{n}\sum_{j=1}^{(m+1)n}N(t_i,t_j)[\mu(t_j)-\mu(t_i)] =-\phi(t_i), 
\quad i=1,2,\ldots,(m+1)n.
\end{equation}
Since $N(s,t)$ is continuous, the term under the summation sign is zero when $j=i$. 
Thus, using the notations $\bx=\mu(\bt)$ and $\by=\phi(\bt)$, the linear system can 
be written as
\begin{equation}\label{e:ie-sys2}
\left(2+\sum_{\begin{subarray}{c} j=1\\j\ne i\end{subarray}}^{(m+1)n}\frac{2\pi}{n}N(t_i,t_j)\right)\bx_i
-\sum_{\begin{subarray}{c} j=1\\j\ne i\end{subarray}}^{(m+1)n}\frac{2\pi}{n}N(t_i,t_j)\bx_j
=-\by_i, \quad i=1,2,\ldots,(m+1)n.
\end{equation}
Let $B$ be an $(m+1)n\times(m+1)n$ matrix with the elements
\begin{equation}\label{e:B}
(B)_{ij} = \left\{
\begin{array}{l@{\hspace{0.5cm}}l}
 \displaystyle 0,  &\text{if $i=j$}, \\[6pt] 
 \displaystyle\frac{2\pi}{n}N(t_i,t_j),     &\text{if $i\ne j$}, 
\quad i,j=1,2,\ldots,(m+1)n.
\end{array}%
\right.
\end{equation}
Thus the $(m+1)n\times(m+1)n$ linear system~(\ref{e:ie-sys2}) 
can be written as
\begin{equation}\label{e:sys}
(2I+\diag(B\bone)-B)\bx=-\by.
\end{equation}

\subsection{Computing the vector $\by$}

In this subsection, we shall present a method for computing the values of the function $\phi(t)$ defined by~(\ref{e:ie-phi}) at the points $t_i$ for $i=1,2,\ldots,(m+1)n$. The method can be used for all H\"older continuous functions $\gamma$ without any differentiability requirement. Thus, the method presented here improves the method presented in~\cite{Nas-siam13} where the function $\gamma$ was assumed to be continuously differentiable. 

We rewrite the index $i$ for $i=1,2,\ldots,(m+1)n$ as 
\[
i=kn+p
\]
where $k=0,1,\ldots,m$ and $p=1,2,\ldots,n$. Hence, by the definitions of the points $t_i$, we need to compute the values 
\begin{equation}\label{e:yi-phik}
y_i = \phi(t_{i})= \phi(t_{kn+p})=\phi_k(s_{k,p}).
\end{equation}
By~(\ref{e:ie-phi}), we have
\[
\phi_k(s_{k,p}) =\int_{J}M(s,t)[\gamma(t)-\gamma_k(s_{kp})]dt
=\displaystyle\sum_{l=0}^{m}\int_{J_l}M(s_{k,p},t)[\gamma_l(t)-\gamma_k(s_{k,p})]dt,
\]
which, in view of~(\ref{e:M-M1}), implies that
\begin{equation}\label{e:phi-kp-1}
\begin{array}{rcl}
\displaystyle\phi_k(s_{k,p}) &= 
&\displaystyle\int_{J_k}\frac{-1}{2\pi}\cot\frac{s_{k,p}-t}{2}[\gamma_k(t)-\gamma_k(s_{k,p})]dt\\
                  & &+\displaystyle\int_{J_k}M_1(s_{k,p},t)[\gamma_k(t)-\gamma_k(s_{k,p})]dt\\
							    & &+\displaystyle\sum_{\begin{subarray}{c} l=0\\l\ne k\end{subarray}}^{m}\int_{J_l}M(s_{k,p},t)
							    [\gamma_l(t)-\gamma_k(s_{k,p})]dt.
\end{array}
\end{equation}
The integral with the cotangent kernel in~(\ref{e:phi-kp-1}) can be discretized by Wittich's method~\cite{Weg05}. The kernel $M_1$ is continuous and the kernel $M$ is continuous for $l\ne k$. So the integrals with the kernels $M_1$ and $M$ 
in~(\ref{e:phi-kp-1}) are discretized by the trapezoidal rule. Hence, we obtain
\begin{equation}\label{e:phi-kp}
\begin{array}{rcl}
  \displaystyle \phi_k(s_{kp})  
	&=&\displaystyle\sum_{q=1}^{n}\left[-(K)_{pq}\right][\gamma_k(s_{kq})-\gamma_k(s_{kp})]  \\[15pt]
	&+&\displaystyle\sum_{q=1}^{n}\frac{2\pi}{n} M_1(s_{kp},s_{kq})[\gamma_k(s_{kq})-\gamma_k(s_{kp})]  \\[15pt]
  & \displaystyle +   &\displaystyle\sum_{\begin{subarray}{c} l=0\\l\ne k\end{subarray}}^{m}\sum_{q=1}^{n}\frac{2\pi}{n}
							M(s_{kp},s_{lq})[\gamma_l(s_{lq})-\gamma_k(s_{kp})],\\
   \end{array}%
\end{equation}
where the elements $(K)_{pq}$ of Wittich's matrix are given by
\[
(K)_{pq} = \left\{
\begin{array}{l@{\hspace{0.5cm}}l}
   0,     &\text{if $p-q$ even}, \\[6pt]
   \frac{2}{n}\cot\frac{(p-q)\pi}{n},     &\text{if $p-q$ odd},  \quad p,q=1,2,\ldots, n.\\
   \end{array}%
   \right. 
\]
Since $(K)_{pq}=0$ where $q=p$ and $M_1(s,t)$ is continuous, the term under the 
first two summation signs in~(\ref{e:phi-kp}) is zero when $q=p$. 
By~(\ref{e:M-M1}), we have for $p\ne q$,
\[
	\frac{2\pi}{n} M_1(s_{kp},s_{kq})
	=\frac{1}{n}\cot\frac{s_{kp}-s_{kq}}{2}+\frac{2\pi}{n} M(s_{kp},s_{kq})=\frac{1}{n} \cot\frac{(p-q)\pi}{n}+\frac{2\pi}{n} M(s_{kp},s_{kq}). 
\]	
Hence, we have for $p\ne q$,
\begin{eqnarray*}
	-(K)_{pq}+\frac{2\pi}{n} M_1(s_{kp},s_{kq})
	&=&-(K)_{pq}+\frac{1}{n} \cot\frac{(p-q)\pi}{n}+\frac{2\pi}{n} M(s_{kp},s_{kq})\\
  &=&(-1)^{p-q}\frac{1}{n}\cot\frac{(p-q)\pi}{n} +\frac{2\pi}{n} M(s_{kp},s_{kq}).
\end{eqnarray*}
Let $L$ be the $n\times n$ matrix whose elements are given by
\[
(L)_{pq} = \left\{
\begin{array}{l@{\hspace{0.5cm}}l}
 \displaystyle 0,  &\text{if $p=q$}, \\[6pt] 
 \displaystyle (-1)^{p-q}\frac{1}{n}\cot\frac{(p-q)\pi}{n},     &\text{if $p\ne q$},
\quad p,q=1,2,\ldots,n. \\
\end{array}%
\right.
\]
Thus we have
\begin{equation}\label{e:K-L}
	-(K)_{pq}+\frac{2\pi}{n} M_1(s_{kp},s_{kq})
	= (L)_{pq}+\frac{2\pi}{n} M(s_{kp},s_{kq}), \quad p\ne q.
\end{equation}

In view of~(\ref{e:K-L}), Eq.~(\ref{e:phi-kp}) can be written as
\begin{equation}\label{e:phi-kp2}
\begin{array}{rcl}
  \displaystyle \phi_k(s_{kp})  
	& = &\displaystyle\sum_{q=1}^{n}(L)_{pq}[\gamma_k(s_{kq})-\gamma_k(s_{kp})]   \\[15pt]
	& + &\displaystyle\sum_{\begin{subarray}{c} q=1\\q\ne p\end{subarray}}^{n}\frac{2\pi}{n} 
	M(s_{kp},s_{kq})[\gamma_k(s_{kq})-\gamma_k(s_{kp})]  \\[15pt]
  & + &\displaystyle\sum_{\begin{subarray}{c} l=0\\l\ne k\end{subarray}}^{m}\sum_{q=1}^{n}\frac{2\pi}{n}
							M(s_{kp},s_{lq})[\gamma_l(s_{lq})-\gamma_k(s_{kp})].\\
   \end{array}%
\end{equation}
Let $D$ be an $(m+1)n\times(m+1)n$ matrix with the elements
\[
(D)_{ij} = \left\{
\begin{array}{l@{\hspace{0.5cm}}l}
 \displaystyle 0,  &\text{if $i=j$}, \\[6pt] 
 \displaystyle\frac{2\pi}{n}M(t_i,t_j),     &\text{if $i\ne j$}, \quad i,j=1,2,\ldots,(m+1)n,\\
\end{array}%
\right.
\]
and $\hat L$ be the $(m+1)n\times(m+1)n$ matrix 
\[
\hat L=\left(\begin{array}{cccc}
L      &O      &\cdots &O \\
O      &L      &\cdots &O \\
\vdots &\vdots &\ddots &\vdots \\
O      &O      &\cdots &L\\
\end{array}\right).
\]
Then, in view of~(\ref{e:ti-skp}), Eq.~(\ref{e:phi-kp2}) can be written as
\begin{equation}\label{e:phi-kp3}
\phi(t_i) =\sum_{j=1}^{(m+1)n}(\tilde L)_{ij}[\gamma(t_j)-\gamma(t_i)] 
          +\sum_{j=1}^{(m+1)n}(D)_{ij}[\gamma(t_{j})-\gamma(t_i)],	
					\quad i=1,2,\ldots,(m+1)n.
\end{equation}
Hence, it follows from~(\ref{e:yi-phik}) and~(\ref{e:phi-kp3}) that the 
$(m+1)n\times1$ vector $\by=\phi(\bt)$ can be written in matrix form as
\begin{equation}\label{e:yo}
\by=D\gamma(\bt)-\diag(D\bone)\gamma(\bt)+\hat L\gamma(\bt)-\diag(\hat L\bone)\gamma(\bt).
\end{equation}

For the vector-matrix product $\hat L\gamma(\bt)$, we have
\begin{equation}\label{e:hL-bt}
\hat L\gamma(\bt) = \left(\begin{array}{cccc}
L      &O      &\cdots &O \\
O      &L      &\cdots &O \\
\vdots &\vdots &\ddots &\vdots \\
O      &O      &\cdots &L\\
\end{array}\right)\gamma(\bt)
= \left(\begin{array}{c}
L\gamma(\bs_1) \\
L\gamma(\bs_2) \\
\vdots \\
L\gamma(\bs_m)\\
\end{array}\right). 
\end{equation} 
The matrix $L$ is circulant since it can be written as
\[
L = \left(\begin{array}{ccccc}
b_1    &b_n    &\cdots &b_3    &b_2 \\
b_2    &b_1    &\ddots &b_4    &b_3 \\
\vdots &\ddots &\ddots &\ddots &\vdots \\
b_{n-1}&b_{n-2}&\ddots &b_1    &b_n\\
b_{n}  &b_{n-1}&\cdots &b_2    &b_1\\
\end{array}\right), 
\]
where $b_1=0$ and 
\[
b_i=(-1)^{i-1}\frac{1}{n}\cot\frac{(i-1)\pi}{n},\quad\text{for}\;\;i=2,3,\ldots,n. 
\]
Thus, the matrix-vector production $L\gamma(\bs_k)$ can be computed in  $O(n\ln n)$ operations using FFT. Using the MATLAB function \verb|fft| for forward FFT and the MATLAB function \verb|ifft| for inverse FFT, the vector $L\gamma(\bs_k)$ is computed by~\cite[p.~92]{Che05}
\begin{equation}\label{e:L-bs}
L\gamma(\bs_k)=\verb|ifft|\left(\verb|fft|(\bvb).*\verb|fft|(\gamma(\bs_k))\right).
\end{equation} 
Hence, the matrix-vector product $\hat L\gamma(\bt)$ can be computed in  
$O((m+1)n\ln n)$ operations. 

Since the \verb|fft| of a constant function is zero, hence, in view of~(\ref{e:L-bs}), we have 
\[
\hat L\bone=0.
\] 
Thus, the vector $\by$ can be written as
\begin{equation}\label{e:y}
\by=D\gamma(\bt)-\diag(D\bone)\gamma(\bt)+\hat L\gamma(\bt).
\end{equation}

For the vector-matrix product $D\gamma(\bt)$, we have for $i=1,2,\ldots,(m+1)n$, 
\begin{eqnarray*}
\sum_{j=1}^{(m+1)n}(D)_{ij}\gamma(t_j) 
          &=&\sum_{\begin{subarray}{c} j=1\\j\ne i\end{subarray}}^{(m+1)n} \frac{2}{n}\Re\left[
          \frac{A(t_i)}{A(t_j)}\frac{\eta'(t_j)}{\eta(t_j)-\eta(t_i)}\right]\gamma(t_j)\\
          &=& -\frac{2}{n}\Re\left[A(t_i)\sum_{\begin{subarray}{c} j=1\\j\ne i\end{subarray}}^{\hat n} 
          \frac{1}{\eta(t_i)-\eta(t_j)}\frac{\eta'(t_j)}{A(t_j)}\gamma(t_j)\right]. 
\end{eqnarray*}
Hence, the matrix-vector product $D\gamma(\bt)$ can be written in terms of
the matrix $E$ as
\begin{equation}\label{e:hD-bt}
D\gamma(\bt) =  -\frac{2}{n} 
\Re\left[A(\bt)\left(E\left(\frac{\eta'(\bt)}{A(\bt)}\gamma(\bt)\right)\right)\right].
\end{equation}
It is clear from~(\ref{e:hD-bt}) that computing $D\gamma(\bt)$ requires one multiplication 
of the matrix $E$ by a vector which can be computed as in~(\ref{e:E-fmm}) by 
FMM in $O((m+1)n)$ operations.  The matrix-vector product $D\bone$ can be also computed in  
$O((m+1)n)$ operations.

Hence, the vector $\by$ in the right-hand side of the linear system~(\ref{e:sys}) can be computed through~(\ref{e:y}) in $O((m+1)n\ln n)$ operations.

\subsection{Multiplication by the coefficient matrix: $2I+\diag(B\bone)-B$}

For multiplying the matrix $B$ by the vector $\bx$, we have for $i=1,2,\ldots,(m+1)n$,
\begin{eqnarray*}
\sum_{j=1}^{(m+1)n}(B)_{ij}\bx_j 
          &=&\sum_{\begin{subarray}{c} j=1\\j\ne i\end{subarray}}^{(m+1)n} \frac{2}{n}\Im\left[
          \frac{A(t_i)}{A(t_j)}\frac{\eta'(t_j)}{\eta(t_j)-\eta(t_i)}\right]\bx_j\\
          &=& -\frac{2}{n}\Im\left[A(t_i)\sum_{\begin{subarray}{c} j=1\\j\ne i\end{subarray}}^{(m+1)n} 
          \frac{1}{\eta(t_i)-\eta(t_j)}\frac{\eta'(t_j)}{A(t_j)}\bx_j\right]. 
\end{eqnarray*}
Hence, the matrix-vector product $B\bx$ can be written in terms of the matrix $E$ as
\begin{equation}\label{e:(I-B)x}
B\bx =-\frac{2}{n}\Im\left[A(\bt)\left(E\left(\frac{\eta'(\bt)}{A(\bt)}\bx\right)\right)\right].
\end{equation}

It is clear from~(\ref{e:(I-B)x}) that multiplying the matrix $B$ by the vector $\bx$ requires one multiplication of the matrix $E$ by a vector which can be computed in $O((m+1)n)$ operations. The matrix-vector product $B\bone$ can be also computed in $O((m+1)n)$ operations. The multiplication of the diagonal matrix $2I+\diag(B\bone)$ by the vector $\bx$ can be computed in $O((m+1)n)$ operations. Thus, the multiplication of the coefficient matrix of the linear system~(\ref{e:sys}), $2I+\diag(B\bone)-B$, by the vector $\bx$ can be computed in $O((m+1)n)$ operations.

\subsection{Solving the linear system}

We use the MATLAB function {\tt gmres} to solve the linear system~(\ref{e:sys}) which can be used with the matrix-vector product function $f_B(\bx)$ defined by
\begin{equation}\label{e:f_B}
f_B(\bx)=(2I+\diag(B\bone)-B)\bx.
\end{equation}
Based on~(\ref{e:(I-B)x}) and~(\ref{e:E-fmm}), the values of the function $f_B(\bx)$ can be computed using the MATLAB function {\tt zfmm2dpart}. The linear system~(\ref{e:sys}) can then be solved using {\tt gmres},
\begin{equation}\label{e:ie-gmres}
\bx = {\tt gmres}(@(\bx)f_B(\bx),-\by,{\tt restart},{\tt tol},{\tt maxit}),
\end{equation}
which restarts every \verb|restart| inner iterations where \verb|tol| is tolerance of the method and \verb|maxit| is the maximum number of outer iterations. By obtaining $\bx$, we obtain an approximation to the solution $\mu$ of the integral equation~(\ref{e:ie}) at the points $\bt$ by $\mu(\bt)=\bx$. 

Since computing the values of the function $f_B(\bx)$ in~(\ref{e:f_B}) requires order $O((m+1)n)$ operations and computing the vector $\by$ requires $O((m+1)n\ln n)$ operations, hence solving the linear system~(\ref{e:sys}) by~(\ref{e:ie-gmres}) requires $O((m+1)n\ln n)$ operations.

\subsection{Computing the piecewise constant function $h$}

In view of~(\ref{e:sys}) and~(\ref{e:y}), the discretizing matrices of the operators $\bI-\bN$ and $\bM$ are $(2I+\diag(B\bone)-B)$ and $D-\diag(D\bone)+\hat L$, respectively. In view of~(\ref{e:h}), the $(m+1)n\times1$ vector $h(\bt)$ which components are the values of $h(t)$ at the points $\bt$ can be approximated by
\begin{equation}\label{e:bh}
h(\bt) = \frac{[D-\diag(D\bone)+\hat L]\mu(\bt)-
[2+\diag(B\bone)-B]\gamma(\bt)}{2}.
\end{equation}
Thus, $h(\bt)$ can be computed in  $O((m+1)n\ln n)$ operations.

Computing the vector $h(\bt)$ in~\cite[Eq.~(62)]{Nas-siam13} required computing $\mu'(\bt)$, i.e., the derivative of the solution of the integral equation~(\ref{e:ie}) at the points $\bt$. Thus, the method presented in~(\ref{e:bh}) improves the method presented in~\cite[Eq.~(62)]{Nas-siam13} since no differentiability of the function $\mu$ is required in~(\ref{e:bh}). 

A MATLAB function \verb|FBIE| for fast solution of the integral equation~(\ref{e:ie}) and fast computing of the function $h$ in~(\ref{e:h}) using the method presented in this section is shown in Figure~\ref{f:FBIE}.

\begin{figure}[!ht]%
\begin{verbatim}
function [mu,h]  =  FBIE(et,etp,A,gam,n,iprec,restart,gmrestol,maxit)
%The function 
%        [mu,h]  =  FBIE(et,etp,A,gam,n,iprec,restart,gmrestol,maxit)
%return the unique solution mu of the integral equation 
%               (I-N)mu=-Mgam 
%and the function 
%                h=[(I-N)gam-Mmu]/2,
%where et is the parameterization of the boundary, etp=et', 
%A=exp(-i\thet)(et-alp) for bounded G and by A=exp(-i\thet) for unbounded
%G, gam is a given function, n is the number of nodes in  each boundary
%component, iprec is the FMM precision flag, restart is the maximum number 
%of GMRES method inner iterations, gmrestol is the tolerance of the GMRES   
%method, and maxit is the maximum number of GMRES method outer iterations
a        = [real(et.') ; imag(et.')];
m        =  length(et)/n-1;
b1       = [etp./A].';
[Ub1]    = zfmm2dpart(iprec,(m+1)*n,a,b1,1);
Eone     = (Ub1.pot).';
b(1,1) = 0;
for k=2:n
    b(k,1) = (-1)^(k+1)*(1/n)*cot(pi*(k-1)/n);
end
mu     = gmres(@(x)fB(x),-fC(gam),restart,gmrestol,maxit);
h       = (fC(mu)-fB(gam))./2;
%%
function  hx  = fB (x)
    bx2   = [x.*etp./A].';
    [Ubx2]= zfmm2dpart(iprec,(m+1)*n,a,bx2,1);
    Ex    = (Ubx2.pot).';
    hx    =  2.*x-(2/n).*imag(A.*Eone).*x+(2/n).*imag(A.*Ex);
end
function  hx = fC (x)
    bx    = [x.*etp./A].';
    [Ubx] = zfmm2dpart(iprec,(m+1)*n,a,bx,1);
    Ex    = (Ubx.pot).';
    for k=1:m+1
        hLx(1+(k-1)*n:k*n,1) = ifft(fft(b).*fft(x(1+(k-1)*n:k*n,1)));
    end
    hx    = -(2/n).*real(A.*Ex)+(2/n).*real(A.*Eone).*x+hLx;    
end
end
\end{verbatim}
\vspace{-0.5cm}
\caption{The MATLAB function {\tt FBIE}.}%
\label{f:FBIE}%
\end{figure}

\section{Solving the integral equation with the adjoint generalized 
Neumann kernel}


\subsection{The integral equation}

The integral equation~(\ref{e:ie*}) will be rewritten using the singularity subtraction in a more suitable form for using the FMM. As for the integral equation~(\ref{e:ie}), the singularity subtraction is useful for solving the integral equation~(\ref{e:ie*}) for domains with corners and for domains with close boundaries. We have
\begin{eqnarray*}
N^\ast(s,t) &=& -\tilde N(s,t) \\
            &=& -\frac{1}{\pi}\Im\left(\frac{\tilde A(s)}{\tilde A(t)}\frac{\eta'(t)}{\eta(t)-\eta(s)}\right) \\
            &=& -\frac{1}{\pi}\Im\left(\frac{\eta'(t)}{\eta(t)-\eta(s)}\right)  
						    +\frac{1}{\pi}\Im\left(\frac{\tilde A(t)-\tilde A(s)}{\eta(t)-\eta(s)}\frac{\eta'(t)}{\tilde A(t)}\right). 
\end{eqnarray*}
Hence,
\begin{equation}\label{e:N*-Nk-Ng}
N^\ast(s,t) = -N_k(s,t)+N_g(s,t)
\end{equation}
where
\begin{equation}\label{e:Nk}
N_k(s,t)=\frac{1}{\pi}\Im\left(\frac{\eta'(t)}{\eta(t)-\eta(s)}\right)
\end{equation}
is the well-known Neumann kernel and 
\begin{equation}\label{e:Ng}
N_g(s,t) 
=\frac{1}{\pi}\Im\left(\frac{\tilde A(t)-\tilde A(s)}{\eta(t)-\eta(s)}\frac{\eta'(t)}{\tilde A(t)}\right)=\frac{1}{\pi}\Im\left(\frac{A(s)\eta'(t)-A(t)\eta'(s)}{A(s)(\eta(t)-\eta(s))}\right). 
\end{equation}
The kernel $N_g(s,t)$ is continuous with 
\begin{equation}\label{e:Ng-tt}
N_g(t,t)=\frac{1}{\pi}\Im\left(\frac{\eta''(t)}{\eta'(t)}-\frac{A'(t)}{A(t)}\right).
\end{equation}

The constant function is an eigenfunction of the Neumann kernel $N_k$ corresponding to the eigenvalue $\lambda=1$ for bounded $G$ and to the eigenvalue $\lambda=-1$ for unbounded $G$~\cite{Nas-amc11,Nas-cr} (see also~\cite{Hen3}). Hence,  
\begin{equation}\label{e:Nk-const}
\int_{J}N_k(s,t)dt=c
\end{equation}
where the constant $c$ is defined by
\begin{equation}\label{e:c}
c= \left\{ \begin{array}{l@{\hspace{0.5cm}}l}
1, &\text{if $G$ is bounded},\\
-1,&\text{if $G$ is unbounded}.
\end{array}
\right.
\end{equation}

Thus, 
\begin{equation}\label{e:N*-const}
\int_{J}N^\ast(s,t)dt=-c+r(s)
\end{equation}
where
\begin{equation}\label{e:r(s)}
r(s)=\int_{J}N_g(s,t)dt.
\end{equation}
Hence, the integral equation~(\ref{e:ie*}) can be written as 
\begin{equation}\label{e:ie*-m}
[1-c+r(s)]\mu(s)+\displaystyle\int_{J}N^\ast(s,t)[\mu(t)-\mu(s)]dt
+\int_{J}\delta(s,t)\mu(t)dt = \gamma(s). 
\end{equation}
The integral equation~(\ref{e:ie*-m}) is valid even if the boundary 
$\Gamma$ is piecewise smooth (see~\cite{Nas-cr}).

\subsection{The Nystr\"om method}

Discretizing the integral in~(\ref{e:ie*-m}) by the 
trapezoidal rule~(\ref{e:trap}) and substituting 
$s=t_i$, we obtain the linear system
\begin{equation}\label{e:ie*-sys}
\begin{array}{rcl}	
\displaystyle [1-c+r(t_i)]\mu(t_i)&+&\displaystyle\frac{2\pi}{n}\sum_{j=1}^{(m+1)n}N(t_j,t_i)[\mu(t_j)-\mu(t_i)] \\
&+&\displaystyle\frac{2\pi}{n}\sum_{j=1}^{(m+1)n}\delta(t_j,t_i)\mu(t_j) =\gamma(t_i), 
\end{array}
\end{equation}
for $i=1,2,\ldots,(m+1)n$. Since $N(s,t)$ is continuous, the term under the summation sign is zero when $j=i$. Thus, the linear system~(\ref{e:ie*-sys}) can be written as
\begin{equation}\label{e:ie*-sys2}
\begin{array}{r}
\displaystyle\left(1-c+r(t_i)-\sum_{\begin{subarray}{c} j=1\\j\ne i\end{subarray}}^{(m+1)n}\frac{2\pi}{n}N(t_j,t_i)\right)\mu(t_i)
+\displaystyle\sum_{\begin{subarray}{c} j=1\\j\ne i\end{subarray}}^{(m+1)n}\frac{2\pi}{n}N(t_j,t_i)\mu(t_j)\\
+\displaystyle\frac{2\pi}{n}\sum_{j=1}^{(m+1)n}\delta(t_j,t_i)\mu(t_j)=\gamma(t_i),
\end{array}
\end{equation}
for $i=1,2,\ldots,(m+1)n$.
Let $P$ be the $n\times n$ matrix with the elements
\[
(P)_{pq} = \frac{1}{n}, \quad p,q=1,2,\ldots,n,
\] 
$\hat P$ be the $(m+1)n\times(m+1)n$ matrix
\[
\hat P=\left(\begin{array}{cccc}
P      &O      &\cdots &O \\
O      &P      &\cdots &O \\
\vdots &\vdots &\ddots &\vdots \\
O      &O      &\cdots &P\\
\end{array}\right),
\]
and $\bve$ be the $(m+1)n\times1$ vector with the elements
\begin{equation}\label{e:ei}
\bve_i=1-c+r(t_i)-\sum_{\begin{subarray}{c} j=1\\j\ne i\end{subarray}}^{(m+1)n}
\frac{2\pi}{n}N(t_j,t_i), \quad i=1,2,\ldots,(m+1)n.
\end{equation}
Let $\bx=\mu(\bt)$ and $B$ be the matrix defined by~(\ref{e:B}). Then, the $(m+1)n\times(m+1)n$ linear system~(\ref{e:ie*-sys2}) can be written as
\begin{equation}\label{e:sys*}
(\diag(\bve)+B^T+\hat P)\bx=\gamma(\bt).
\end{equation} 

\subsection{Computing the vector $\bve$}

By approximating the values of the function $r$ at the points $t_i$, i.e., 
\[
r(t_i)=\int_{J}N_g(t_i,t)dt, \quad i=0,1,\ldots,(m+1)n,
\]
by the trapezoidal rule~(\ref{e:trap}), we obtain
\begin{equation}\label{e:ri}
r(t_i) =\displaystyle\frac{2\pi}{n}\sum_{j=1}^{(m+1)n}N_g(t_i,t_j) 
=\frac{2\pi}{n}\sum_{\begin{subarray}{c} j=1\\j\ne i\end{subarray}}^{(m+1)n}N_g(t_i,t_j)
+\frac{2\pi}{n}N_g(t_i,t_i).
\end{equation}
By~(\ref{e:N*-Nk-Ng}), we have
\begin{equation}\label{e:ri2}
\frac{2\pi}{n}\sum_{\begin{subarray}{c} j=1\\j\ne i\end{subarray}}^{(m+1)n}
N_g(t_i,t_j) =\frac{2\pi}{n}\sum_{\begin{subarray}{c} j=1\\j\ne i\end{subarray}}^{(m+1)n}
N_k(t_i,t_j)
+\frac{2\pi}{n}\sum_{\begin{subarray}{c} j=1\\j\ne i\end{subarray}}^{(m+1)n}
N(t_j,t_i).
\end{equation}
Thus, it follows from~(\ref{e:ei}),~(\ref{e:ri}) and~(\ref{e:ri2}) that
\begin{equation}\label{e:ei2}
\bve_i=1-c+
\frac{2\pi}{n}\sum_{\begin{subarray}{c} j=1\\j\ne i\end{subarray}}^{(m+1)n}N_k(t_i,t_j) 
+\frac{2\pi}{n}N_g(t_i,t_i), \quad i=0,1,\ldots,(m+1)n.
\end{equation}

For $i=0,1,\ldots,(m+1)n$, by~(\ref{e:Ng-tt}) and~(\ref{e:Nk}), we have
\[
\frac{2\pi}{n}N_g(t_i,t_i)= \frac{2}{n}\Im\left(\frac{\eta''(t_i)}{\eta'(t_i)}-\frac{A'(t_i)}{A(t_i)}\right)
\]
and
\[
\frac{2\pi}{n}\sum_{\begin{subarray}{c} j=1\\j\ne i\end{subarray}}^{(m+1)n}N_k(t_i,t_j)
=\frac{2}{n}\sum_{\begin{subarray}{c} j=1\\j\ne i\end{subarray}}^{(m+1)n}
\Im\left(\frac{\eta'(t_j)}{\eta(t_j)-\eta(t_i)}\right) =-\frac{2}{n}\Im\left[\sum_{\begin{subarray}{c} j=1\\j\ne i\end{subarray}}^{(m+1)n}
\frac{1}{\eta(t_i)-\eta(t_j)}\eta'(t_j)
\right].
\]
Hence, in view of~(\ref{e:fmm-E}) and~(\ref{e:ei2}), 
the vector $\bve$ can be written as
\begin{equation}\label{e:bve}
\bve = 1-c-\frac{2}{n}\Im\left[E\eta'(\bt)\right]
+\frac{2}{n}\Im\left[\frac{\eta''(\bt)}{\eta'(\bt)}-\frac{A'(\bt)}{A(\bt)}\right].
\end{equation}
Thus, computing the vector $\bve$ requires one multiplication of the matrix $E$ by a 
vector which can be computed in $O((m+1)n)$ operations.  

\subsection{Multiplication by the coefficient matrix: $\diag(\bve)+B^T+\hat P$}

For multiplying the matrix $B^T$ by the vector $\bx$, we have for $i=1,2,\ldots,(m+1)n$,
\begin{eqnarray*}
\sum_{j=1}^{(m+1)n}(B^T)_{ij}\bx_j &=&\sum_{j=1}^{(m+1)n}(B)_{ji}\bx_j \\
          &=&\sum_{\begin{subarray}{c} j=1\\j\ne i\end{subarray}}^{(m+1)n} \frac{2}{n}\Im\left[
          \frac{A(t_j)}{A(t_i)}\frac{\eta'(t_i)}{\eta(t_i)-\eta(t_j)}\right]\bx_j\\
          &=& \frac{2}{n}\Im\left[\frac{\eta'(t_i)}{A(t_i)}\sum_{\begin{subarray}{c} j=1\\j\ne i\end{subarray}}^{(m+1)n} 
          \frac{1}{\eta(t_i)-\eta(t_j)}A(t_j)\bx_j\right]. 
\end{eqnarray*}
Hence, the matrix-vector product $B^T\bx$ can be written in terms of the matrix $E$ as
\begin{equation}\label{e:BTx}
B^T\bx =\frac{2}{n}\Im\left[\frac{\eta'(\bt)}{A(\bt)}\left(E\left(A(\bt)\bx\right)\right)\right].
\end{equation}
Thus, multiplying the matrix $B^T$ by the vector $\bx$ requires one multiplication of the matrix $E$ by a vector which can be computed in $O((m+1)n)$ operations.

For the matrix-vector product $\hat P\bx$, we rewrite the vector $\bx$ as 
\[
\bx=\left(\begin{array}{c}
\bx_0 \\
\bx_1 \\
\vdots \\
\bx_m\\
\end{array}\right).
\]
where each of the vectors $\bx_j$, $j=0,1,\ldots,m$, is an $n\times 1$ vector. Then the matrix-vector product $\hat P\bx$ can be computed by
\[
\hat P\bx=\left(\begin{array}{c}
P\bx_0 \\
P\bx_1 \\
\vdots \\
P\bx_m\\
\end{array}\right)=
\left(\begin{array}{c}
\frac{1}{n}\sum_{p=1}^{n}\left(\bx_0\right)_p \\[4pt]
\frac{1}{n}\sum_{p=1}^{n}\left(\bx_1\right)_p \\
\vdots \\
\frac{1}{n}\sum_{p=1}^{n}\left(\bx_m\right)_p\\
\end{array}\right).
\]
Thus the multiplication of the matrix $\hat P$ by the vector $\bx$ can be computed in $O((m+1)n)$ operations. The multiplication of the diagonal matrix $\diag(\bve)$ by the vector $\bx$ can also be computed in $O((m+1)n)$ operations. Hence, the multiplication of the coefficient matrix $\diag(\bve)+B^T+\hat P$ of the linear system~(\ref{e:sys*}) by the vector $\bx$ can be computed in  $O((m+1)n)$ operations.

\subsection{Solving the linear system~(\ref{e:sys*})}

The linear system~(\ref{e:sys*}) will be solved using the MATLAB function {\tt gmres} with the matrix-vector product function 
\begin{equation}\label{e:g_B}
g_B(\bx)=(\diag(\bve)+B^T+\hat P)\bx.
\end{equation}
Based on~(\ref{e:BTx}) and~(\ref{e:E-fmm}), the values of the function $g_B(\bx)$ can be computed using the MATLAB function {\tt zfmm2dpart}. The linear system~(\ref{e:sys*}) can then be solved using {\tt gmres},
\begin{equation}\label{e:ie*-gmres}
\bx = {\tt gmres}(@(\bx)g_B(\bx),\gamma(\bt),{\tt restart},{\tt tol},{\tt maxit}).
\end{equation}
By obtaining $\bx$, we obtain an approximation to the solution $\mu$ of the integral equation~(\ref{e:ie*}) at the points $\bt$ by $\mu(\bt)=\bx$. 

In contrast to the integral equation with the generalized Neumann kernel~(\ref{e:ie}), the right-hand side of the integral equation with the adjoint generalized Neumann kernel~(\ref{e:ie*}) is given explicitly. Since computing the values of the function $g_B(\bx)$ in~(\ref{e:g_B}) requires $O((m+1)n)$ operations, hence solving the linear system~(\ref{e:sys*}) by~(\ref{e:ie*-gmres}) requires $O((m+1)n)$ operations. 

A MATLAB function \verb|FBIEad| for fast solution of the integral equation~(\ref{e:ie*}) using the method presented in this section is shown in Figure~\ref{f:FBIEad}.

\begin{figure}[!ht]%
\begin{verbatim}
function mu=FBIEad(et,etp,etpp,A,Ap,gam,n,c,iprec,restart,gmrestol,maxit)
%The function 
%        mu=FBIEad(et,etp,etpp,A,Ap,gam,n,c,iprec,restart,gmrestol,maxit)
%returns the unique solution mu of the integral equation 
%               (I+N*+J)mu=gam 
%where et is the parameterization of the boundary, etp=et', etpp=et'',
%A=exp(-i\thet)(et-alp) for bounded G and A=exp(-i\thet) for unbounded
%G, gam is a given function, n is the number of nodes in each boundary
%component, c=1 for bounded G and c=-1 for unbounded G, iprec is the
%FMM precision flag, restart is the maximum number of the GMRES method 
%inner iterations, gmrestol is the tolerance of the GMRES method, and
%maxit is the maximum number of GMRES method outer iterations
a      = [real(et.') ; imag(et.')];
m      =  length(et)/n-1;
[Uetp] = zfmm2dpart(iprec,(m+1)*n,a,etp.',1);
Eetp   = (Uetp.pot).';
e      = 1-c-(2/n).*imag(Eetp)+(2/n).*imag(etpp./etp-Ap./A);     
mu     = gmres(@(x)gB(x),gam,restart,gmrestol,maxit);
%%
function  hx  = gB (x)
    for k=1:m+1
        hPx(1+(k-1)*n:k*n,1) = (1/n)*sum(x(1+(k-1)*n:k*n,1));
    end
    [UAx]  = zfmm2dpart(iprec,(m+1)*n,a,[A.*x].',1);
    EAx    = (UAx.pot).';
    Btx    = (2/n).*imag((etp./A).*EAx); 
    hx    =  e.*x+Btx+hPx;
end
end
\end{verbatim}
\vspace{-0.5cm}
\caption{The MATLAB function {\tt FBIEad}.}%
\label{f:FBIEad}%
\end{figure}

\section{Computing the Cauchy integral formula}

The solutions of the boundary integral equations~(\ref{e:ie}) and~(\ref{e:ie*}) provide us with the values of the conformal mapping and the solution of the boundary value problem on the boundary~$\Gamma$. Computing the values of the conformal mapping and the solution of the boundary value problems for interior points $z\in G$ required computing the Cauchy integral formula
\begin{subequations}\label{e:cauchy}
\renewcommand{\theequation}{\theparentequation\alph{equation}}
\begin{equation}\label{e:cauchy-b}
f(z)=\frac{1}{2\pi\i}\int_{\Gamma}\frac{f(\eta)}{\eta-z}d\eta, \quad z\in G,
\end{equation}
for bounded $G$ and
\begin{equation}\label{e:cauchy-u}
f(z)=f(\infty)+\frac{1}{2\pi\i}\int_{\Gamma}\frac{f(\eta)}{\eta-z}d\eta, \quad z\in G,
\end{equation}
\end{subequations}
for unbounded $G$ where $f$ is an analytic function on $G$ with known boundary values on the boundary $\Gamma$. In this section, we shall review an accurate and fast numerical method for computing the Cauchy integral formulas~(\ref{e:cauchy-b}) and~(\ref{e:cauchy-u}) from~\cite{Aus13,Hel-Oja,Nas-siam13}.

The integral in~(\ref{e:cauchy}) can be approximated by the trapezoidal rule~(\ref{e:trap}). However, the integrand in~(\ref{e:cauchy}) becomes nearly singular for points $z\in G$ which are close to the boundary~$\Gamma$. For such case, the singularity subtraction can be used to obtain accurate results~\cite{Hel-Oja}. See also~\cite{Aus13}. The FMM can be used for fast computing of the values of the function $f(z)$~\cite{Nas-siam13}. 

Suppose that $f(z)$ is analytic in a domain $\hat G$ contains $G\cup\Gamma$. Thus, the integrand function in~(\ref{e:cauchy}), i.e.,
\begin{equation}\label{e:cau-1}
\frac{f(\eta)}{\eta-z}
\end{equation}
has a ploe at $\eta=z\in\hat G$. For $z\in G$, we have
\begin{subequations}\label{e:cauchy-1}
\renewcommand{\theequation}{\theparentequation\alph{equation}}
\begin{equation}\label{e:cauchy-1-b}
\frac{1}{2\pi\i}\int_{\Gamma}\frac{1}{\eta-z}d\eta=1, 
\end{equation}
for bounded $G$ and
\begin{equation}\label{e:cauchy-1-u}
\frac{1}{2\pi\i}\int_{\Gamma}\frac{1}{\eta-z}d\eta=0, 
\end{equation}
\end{subequations}
for unbounded $G$. Thus, the Cauchy integral formula~(\ref{e:cauchy}) can be then written for $z\in G$ as
\begin{subequations}\label{e:cauchy2}
\renewcommand{\theequation}{\theparentequation\alph{equation}}
\begin{equation}\label{e:cauchy2-b}
\frac{1}{2\pi\i}\int_{\Gamma}\frac{f(\eta)-f(z)}{\eta-z}d\eta=0, 
\end{equation}
for bounded $G$ and
\begin{equation}\label{e:cauchy2-u}
f(z)=f(\infty)+\frac{1}{2\pi\i}\int_{\Gamma}\frac{f(\eta)-f(z)}{\eta-z}d\eta, 
\end{equation}
\end{subequations}
for unbounded $G$. Thus $\eta=z$ is not a ploe of the integrand function in the new formula~(\ref{e:cauchy2}) since the integrand
\begin{equation}\label{e:cau-2}
\frac{f(\eta)-z}{\eta-z}
\end{equation}
is an analytic function of $\eta\in\hat G$. Hence, accurate results can be obtained if the trapezoidal rule is used to discretize the integral in~(\ref{e:cauchy2}) (ses~\cite{Aus13} for more details). 

By discretizing the integral in~(\ref{e:cauchy2}) using the trapezoidal rule~(\ref{e:trap}), we obtain for $z\in G$,
\begin{subequations}\label{e:cauchy3}
\renewcommand{\theequation}{\theparentequation\alph{equation}}
\begin{equation}\label{e:cauchy3-b}
\frac{2\pi}{n}\frac{1}{2\pi\i}\sum_{j=1}^{(m+1)n}\frac{f(\eta(t_j))-f(z)}{\eta(t_j)-z}\eta'(t_j)
\approx0
\end{equation}
for bounded $G$ and
\begin{equation}\label{e:cauchy3-u}
f(z)\approx f(\infty)+
\frac{2\pi}{n}\frac{1}{2\pi\i}\sum_{j=1}^{(m+1)n}\frac{f(\eta(t_j))-f(z)}{\eta(t_j)-z}\eta'(t_j)
\end{equation}
\end{subequations}
for unbounded $G$. Consequently, the values of the function $f(z)$ is given for $z\in G$ by
\begin{subequations}\label{e:cauchy4}
\renewcommand{\theequation}{\theparentequation\alph{equation}}
\begin{equation}\label{e:cauchy4-b}
f(z)\approx \frac{\displaystyle \sum_{j=1}^{(m+1)n}\frac{f(\eta(t_j))\eta'(t_j)}{\eta(t_j)-z}}
{\displaystyle \sum_{j=1}^{(m+1)n}\frac{\eta'(t_j)}{\eta(t_j)-z}}
\end{equation}
for bounded $G$ and
\begin{equation}\label{e:cauchy4-u}
f(z)\approx \frac{\displaystyle f(\infty)+
\frac{1}{n\i}\sum_{j=1}^{(m+1)n}\frac{f(\eta(t_j))\eta'(t_j)}{\eta(t_j)-z}}
{\displaystyle 1+\frac{1}{n\i}\sum_{j=1}^{(m+1)n}\frac{\eta'(t_j)}{\eta(t_j)-z}}
\end{equation}
\end{subequations}
for unbounded $G$. 

For a given positive integer $\hat n$, let $z_1,z_2,\ldots,z_{\hat n}$ be given points in $G$ and $\bz$ be the $\hat n\times1$ complex vector $\bz=(z_1,z_2,\ldots,z_{\hat n})$. Hence, the values $f(z_i)$ for $i=1,2,\ldots,\hat n$ can be computed from
\begin{subequations}\label{e:cauchy5}
\renewcommand{\theequation}{\theparentequation\alph{equation}}
\begin{equation}\label{e:cauchy5-b}
f(z_i)\approx \frac{\displaystyle\sum_{j=1}^{(m+1)n}\frac{1}{z_i-\eta(t_j)}f(\eta(t_j))\eta'(t_j)}
{\displaystyle\sum_{j=1}^{(m+1)n}\frac{1}{z_i-\eta(t_j)}\eta'(t_j)}
\end{equation}
for bounded $G$ and
\begin{equation}\label{e:cauchy5-u}
f(z_i)\approx \frac{\displaystyle f(\infty)-
\frac{1}{n\i}\sum_{j=1}^{(m+1)n}\frac{1}{z_i-\eta(t_j)}f(\eta(t_j))\eta'(t_j)}
{\displaystyle 1-\frac{1}{n\i}\sum_{j=1}^{(m+1)n}\frac{1}{z_i-\eta(t_j)}\eta'(t_j)}
\end{equation}
\end{subequations}
for unbounded $G$. Using the matrix $F$ defined by~(\ref{e:F}), the summations~(\ref{e:cauchy5}) can be written as a matrix-vector product  
\begin{subequations}\label{e:cauchy6}
\renewcommand{\theequation}{\theparentequation\alph{equation}}
\begin{equation}\label{e:cauchy6-b}
f(\bz)\approx \frac{F[f(\eta(\bt))\eta'(\bt)]}{F[\eta'(\bt)]}
\end{equation}
for bounded $G$ and
\begin{equation}\label{e:cauchy6-u}
f(\bz)\approx \frac{\displaystyle f(\infty)-
\frac{1}{n\i}F[f(\eta(\bt))\eta'(\bt)]}
{\displaystyle 1-\frac{1}{n\i}F[\eta'(\bt)]}
\end{equation}
\end{subequations}
for unbounded $G$. Hence, computing the vector $f(\bz)$ in~(\ref{e:cauchy6}) requires two multiplications of the matrix $F$ by a vector which can be computed as in~(\ref{e:F-fmm}) by FFF in order $O((m+1)n+\hat n)$ operations.

A MATLAB function \verb|FCAU| for fast computing of $f(\bz)$ using the above described method is presented in Figure~\ref{f:CAU}. 

\begin{figure}%
\begin{verbatim}
function  fz  = FCAU (et,etp,f,z,n,finf)
%The function 
%        fz  = FCAU (et,etp,f,z,n,finf)
%returns the values of the analytic function f computed using the Cauchy
%integral formula at interior vector of points z, where et is the
%parameterization of the boundary, etp=et', finf is the values of f at 
%infinity for unbounded G, n is the number of nodes in each boundary
%component. 
vz    = [real(z) ; imag(z)];       % target
nz    = length(z);                 % ntarget
a     = [real(et.') ; imag(et.')]; % source
tn    = length(et);                % nsource=(m+1)n
iprec = 4;                         %- FMM precision flag
bf    = [f.*etp].';
[Uf]  = zfmm2dpart(iprec,tn,a,bf,0,0,0,nz,vz,1,0,0);
b1    = [etp].';
[U1]  = zfmm2dpart(iprec,tn,a,b1,0,0,0,nz,vz,1,0,0);
if( nargin == 4 ) 
    fz    = (Uf.pottarg)./(U1.pottarg);
end
if( nargin == 6 ) 
    fz= (finf-(Uf.pottarg)./(n*i))./(1-(U1.pottarg)./(n*i));
end
end
\end{verbatim}
\vspace{-0.5cm}
\caption{The MATLAB function {\tt FCAU}.}%
\label{f:CAU}%
\end{figure}

\section{Domains with piecewise smooth boundary}

Suppose that $\gamma(t)$ is smooth in each interval $J_j$ except at $p_j\ge1$ 
points 
\[
c_{j,k}=(k-1)\frac{2\pi}{p_j}\in J_j, \quad k=1,2,\dots,p_j,\quad j=0,1,\ldots,m. 
\]
Suppose that $\omega(t)$ is the bijective, strictly monotonically increasing and infinitely differentiable function defined by~\cite{Kre90}
\begin{equation}\label{e:omega}
\omega(t)=2\pi\frac{[v(t)]^p}{[v(t)]^p+[v(2\pi-t)]^p},
\end{equation}
where
\begin{equation}\label{e:omega-v}
v(t)=\left( \frac{1}{p}-\frac{1}{2}\right)\left(\frac{\pi-t}{\pi}\right)^3 
+\frac{1}{p}\frac{t-\pi}{\pi}+\frac{1}{2},
\qquad t\in [0,2\pi].
\end{equation}
The grading parameter $p$ is an integer such that $p\ge2$. 

We define a function $\delta_j(t)$,
\[
\delta_j(t)\;:\;[0,2\pi]\to[0,2\pi], 
\]
by
\begin{equation}\label{e:delta_j}
\delta_j(t)=\frac{1}{p_j}\omega(p_j\,(t-c_{j,k}))+c_{j,k}, \quad 
t\in\left[c_{j,k},c_{j,k+1}\right],
\end{equation}
for $k=1,2,\ldots,p_j$ and $j=0,1,\ldots,m$. Then the function $\delta_j$ satisfies
\[
\delta'_j\left(c_{j,k}\right)=0,\quad k=1,2,\ldots,p_j,
\]
\[
\delta'_j(t)\ne0\quad \text{for all}\;\; t\in J_j-\{c_{j,1},c_{j,2},\ldots,c_{j,p_j}\}.
\]

Thus, to compute the integral $\int_j\gamma(t)dt$, we introduce the substitution $t=\delta(\tau)$ to obtain
\[
\int_J\gamma(t)dt = \int_J\gamma(\delta(\tau))\delta'(\tau)d\tau = \int_J\hat\gamma(\tau)d\tau
\]
where
\[
\hat\gamma(\tau) = \gamma(\delta(\tau))\delta'(\tau).
\]
The function $\hat\gamma$ is smooth on $J$ and satisfies $\hat\gamma(0)=\hat\gamma(2\pi)=0$. 
Hence, applying the trapezoidal rule~(\ref{e:trap}) to the 
transform integral yields
\begin{equation}\label{e:trap-cr}
\int_J\gamma(t)dt \approx \frac{2\pi}{n}\sum_{k=0}^{m}\sum_{p=1}^{n}\hat\gamma(s_{k,p})
=\frac{2\pi}{n}\sum_{j=1}^{(m+1)n}\hat\gamma(t_{j})
=\frac{2\pi}{n}\sum_{j=1}^{(m+1)n}\gamma(\delta(t_{j}))\delta'(t_{j}).
\end{equation}

Now, suppose that each boundary component $\Gamma_j$ contains $p_j$ corner points located at $\eta_j(c_{j,k})$, $k=1,2,\dots,p_j$, $j=0,1,\ldots,m$. Then the integral equations~(\ref{e:ie-m}) and~(\ref{e:ie*-m})  can be solved accurately by discretizing the integrals in the integral equations~(\ref{e:ie-m}) and~(\ref{e:ie*-m}) by the trapezoidal rule~(\ref{e:trap-cr}) (see~\cite{Kre90,Nas-cr,Yun-sp} for more details). 

An equivalent method for discretizing the integrals in the integral equations~(\ref{e:ie-m}) and~(\ref{e:ie*-m}) using the trapezoidal rule~(\ref{e:trap-cr}) is to choose a piecewise smooth parametrization $\zeta_j(t)$ of the boundary component $\Gamma_j$ then defining a parametrization $\eta_j(t)$ of $\Gamma_j$ by
\begin{equation}\label{e:zet-et}
\eta_j(t)=\zeta_j(\delta_j(t)), \quad j=0,1,\ldots,m.
\end{equation}
The integrals in the integral equations~(\ref{e:ie-m}) and~(\ref{e:ie*-m}) can then be discretized accurately by the trapezoidal rule~(\ref{e:trap}). See~\cite{Els-Gra,Yun-sp} for more details.

\section{Numerical Examples}

To test the performance of the presented functions \verb|FBIE| and \verb|FBIEad|, four numerical examples are presented. In the first example, we consider bounded and unbounded multiply connected domains of high connectivity. For the bounded domain, half of the boundaries are piecewise smooth boundaries. In the second, third and fourth examples, we consider a bounded and an unbounded circular domain of connectivity $5$ with variable distance $\varepsilon$ between the boundaries. In the second example, we show the accuracy of the functions \verb|FBIE| and \verb|FBIEad| for domains whose boundaries are very close to each other. In the third and fourth examples, we test the performance of the functions \verb|FBIE| and \verb|FBIEad| by computing the conformal mapping from bounded and unbounded domains with close boundaries onto the unit disc with circular slits and the unit disc with both circular and radial slits. 

\begin{figure}[ht]%
\centerline{\scalebox{0.85}{\includegraphics{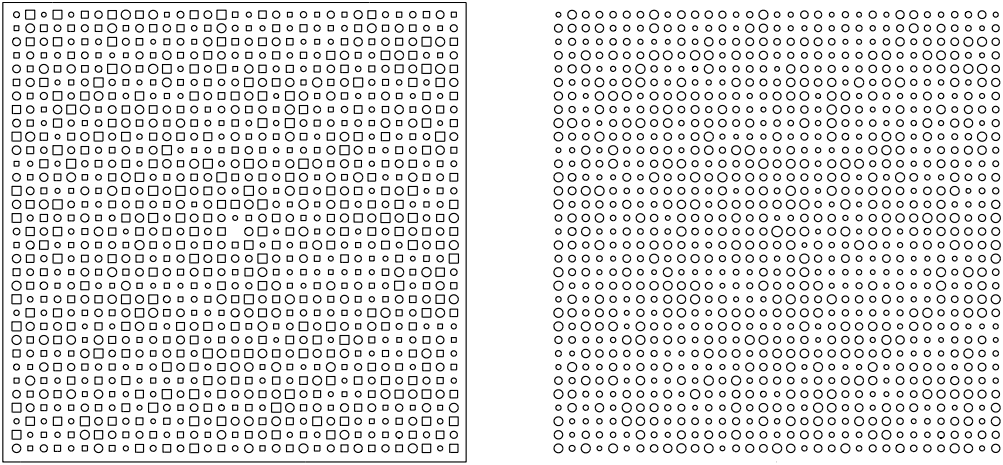}}}
\caption{\rm The domains of Example~1.} 
\label{f:ex1-fig}
\end{figure}

\vspace{0.25cm}
{\noindent\bf Example 1.}
We consider a bounded and an unbounded multiply connected domains $G$ of connectivity $1089$ ($m=1088$) 
(see Fig.~\ref{f:ex1-fig}). The boundary of the bounded domain consists of $544$ circles and $545$ squares including the external boundary. The boundary of the unbounded domain consists of $1089$ circles.

We assume that $f(z)$ is defined for bounded $G$ by
\begin{equation}\label{e:ex-fb}
f(z)=\sin z+\frac{1}{z-2},
\end{equation}
and for unbounded $G$ by
\begin{equation}\label{e:ex-fu}
f(z)=\frac{1}{z}-\sin\frac{1}{z}.
\end{equation}
The function $f(z)$ is an analytic function in $G$ with $f(\infty)=0$ for unbounded $G$. We assume also that
\begin{equation}\label{e:ex-gam}
\gamma(t)=\Re[A(t)f(\eta(t))], 
\end{equation}
where the function $A$ is defined by~(\ref{e:A}) with $\alpha=0$ for bounded $G$ and $\theta(t)=(\theta_0,\theta_1,\ldots,\theta_m)$ where
\[
\theta_j=\frac{2j\pi}{m}, \quad j=0,1,,\ldots,m.
\]

\begin{figure}[ht]%
\centerline{\scalebox{0.7}{\includegraphics{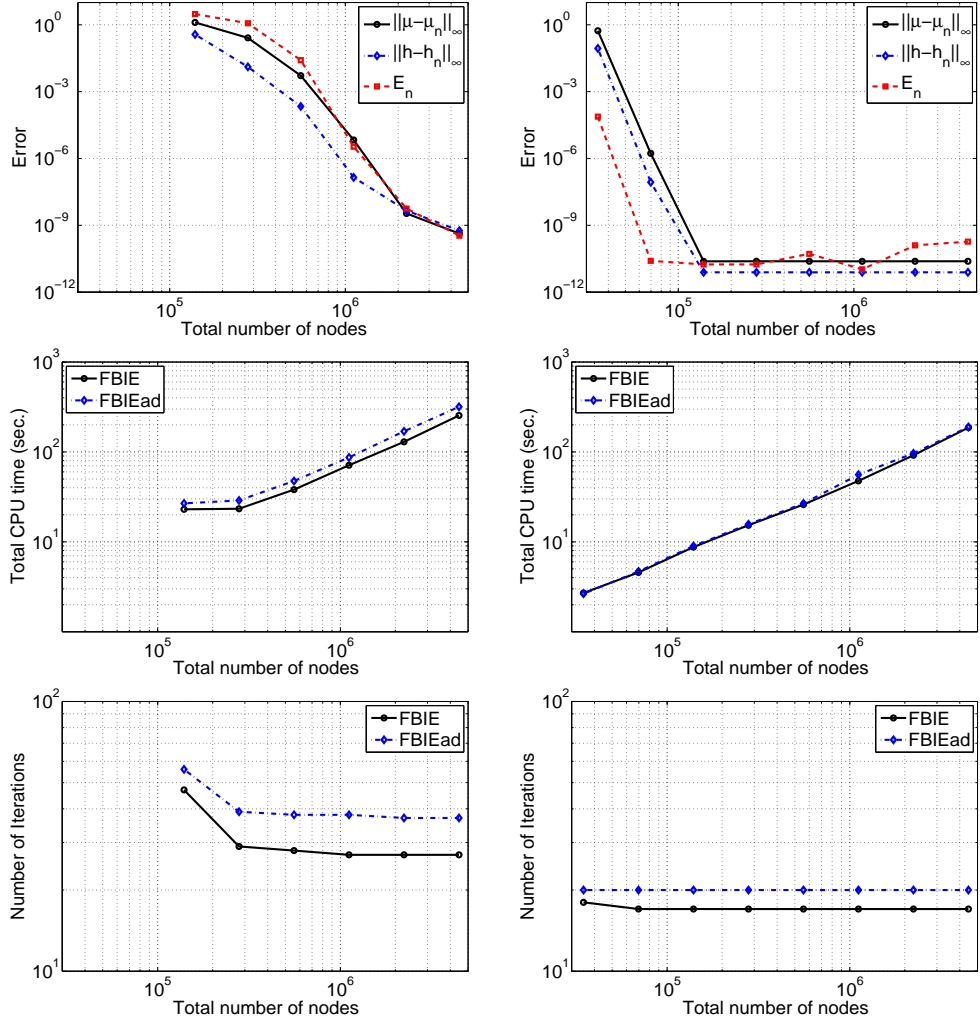}}}
\caption{\rm The maximum error norms $\|\mu-\mu_n\|_\infty$ and $\|h-h_n\|_\infty$, the error $E_n$, the total CPU time, and the number of GMRES iterations vs. the total number of nodes for the bounded domain (left) and the unbounded domain (right).} 
\label{f:ex1-err}
\end{figure}

We consider the numerical solution of the uniquely solvable integral equation with 
the generalized Neumann kernel
\begin{equation}\label{e:ex-ie}
(\bI+\bN)\mu=-\bM\gamma
\end{equation}
and the numerical computing the function $h$ given by
\begin{equation}\label{e:ex-h}
h = [\bM\mu-(\bI-\bN)\gamma]/2.
\end{equation}
The exact solution of the integral equation~(\ref{e:ex-ie}) is 
\[
\mu(t)=\Im[A(t)f(\eta(t))], \quad t\in J,
\]
and the exact value of the function $h$ in~(\ref{e:ex-h}) is~\cite{Weg-Nas,Nas-cmft09} 
\[
h(t)=0, \quad t\in J,
\]
Suppose that $\mu_n$ and $h_n$ are the approximate solutions obtained using the function \verb|FBIE|. The values of the maximum error norms $\|\mu-\mu_n\|_\infty$ and $\|h-h_n\|_\infty$ vs. the total number of nodes are shown in Fig.~\ref{f:ex1-err}. 

We shall also consider the numerical solution of the uniquely solvable integral equation with the adjoint generalized Neumann kernel
\begin{equation}\label{e:ex-ie*}
(\bI+\bN^\ast+\bJ)\phi=1.
\end{equation}
We do not have the exact solution of~(\ref{e:ex-ie*}). However, its unique solution $\phi$ satisfies~\cite{Nas-amc11}
\[
(\bI+\bN^\ast)\phi=0, \quad \bJ\phi=1.
\]
Since the Riemann-Hilbert problem 
\[
\Re[A(t)f(\eta(t))]=\gamma(t)
\]
is solvable for the function $\gamma$ given by~(\ref{e:ex-gam}), hence we have~\cite{Nas-amc11}
\[
\int_J\gamma(t)\phi(t)dt=0.
\]
Let $\phi_{n}$ be an approximation to the unique solution of the integral equation~(\ref{e:ex-ie*}) obtained using the function \verb|FBIEad|. Then,
\[
\int_J\gamma(t)\phi_n(t)dt
\approx \frac{2\pi}{n}\sum_{j=1}^{(m+1)n}\gamma(t_j)\phi_n(t_j)
\approx 0.
\]
We define the error in $\phi_n$ by
\[
E_n=\left|\frac{2\pi}{n}\sum_{j=1}^{(m+1)n}\gamma(t_j)\phi_n(t_j)\right|.
\]
The values of the error $E_n$ vs. the total number of nodes are shown in Fig.~\ref{f:ex1-err}. 

\begin{figure}[ht]%
\centerline{\scalebox{0.7}{\includegraphics{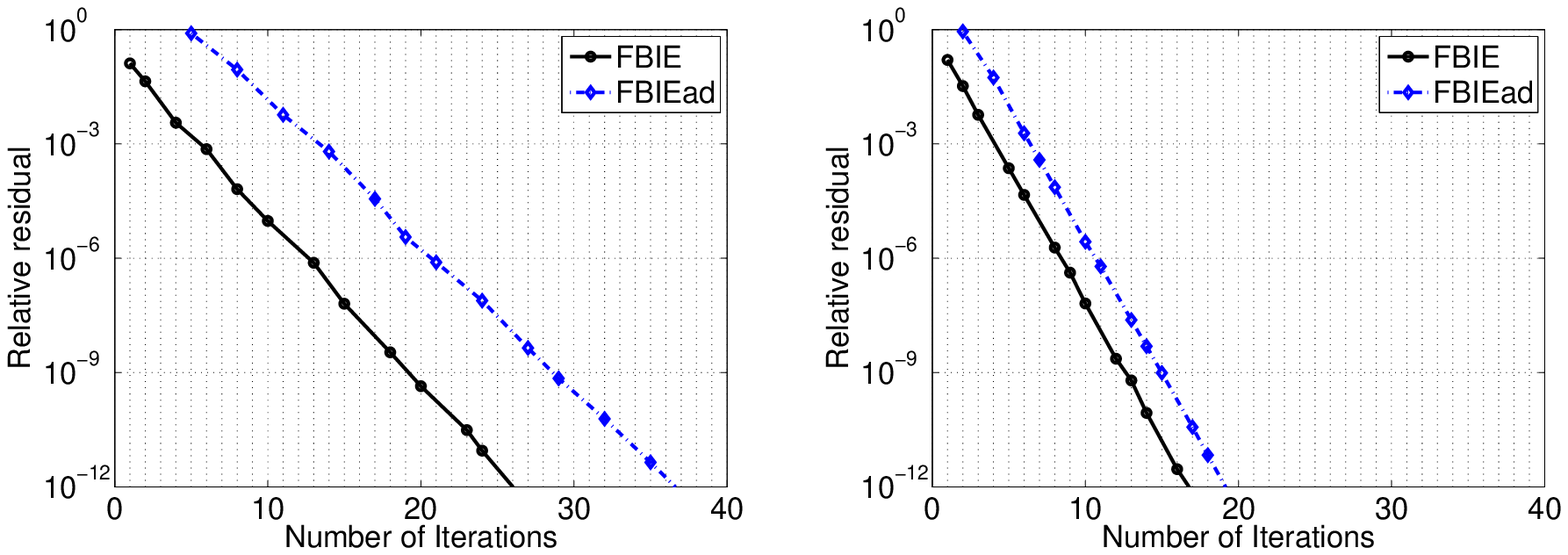}}}
\caption{\rm The relative residual vs. the number of GMRES iterations for the bounded domain (left) and the unbounded domain (right).} 
\label{f:ex1-gmres}
\end{figure}

Figure~\ref{f:ex1-err} shows also the total CPU time (in seconds) and the Number of GMRES iterations required to obtain the approximate solutions $\mu_n$, $h_n$ using the function \verb|FBIE| and $\phi_n$ using the function \verb|FBIEad| vs. the total number of nodes. The numerical results are obtained with ${\tt iprec}=4$, ${\tt restart}=10$, ${\tt gmrestol}=10^{-12}$, and ${\tt maxit}=10$. The relative residual vs. the number of iterations of GMRES obtained with $n=4096$ (total number of nodes is $4460544$) is shown in Figure~\ref{f:ex1-gmres}. It is clear from the Figures~\ref{f:ex1-err} and~\ref{f:ex1-gmres} that the accuracy of the numerical results for the unbounded domain is better than for the bounded domain. This is expected since all the boundaries of the unbounded domain are smooth and half of the boundaries of the bounded domain have corners.

\begin{figure}[ht]%
\centerline{\scalebox{0.3}{\includegraphics{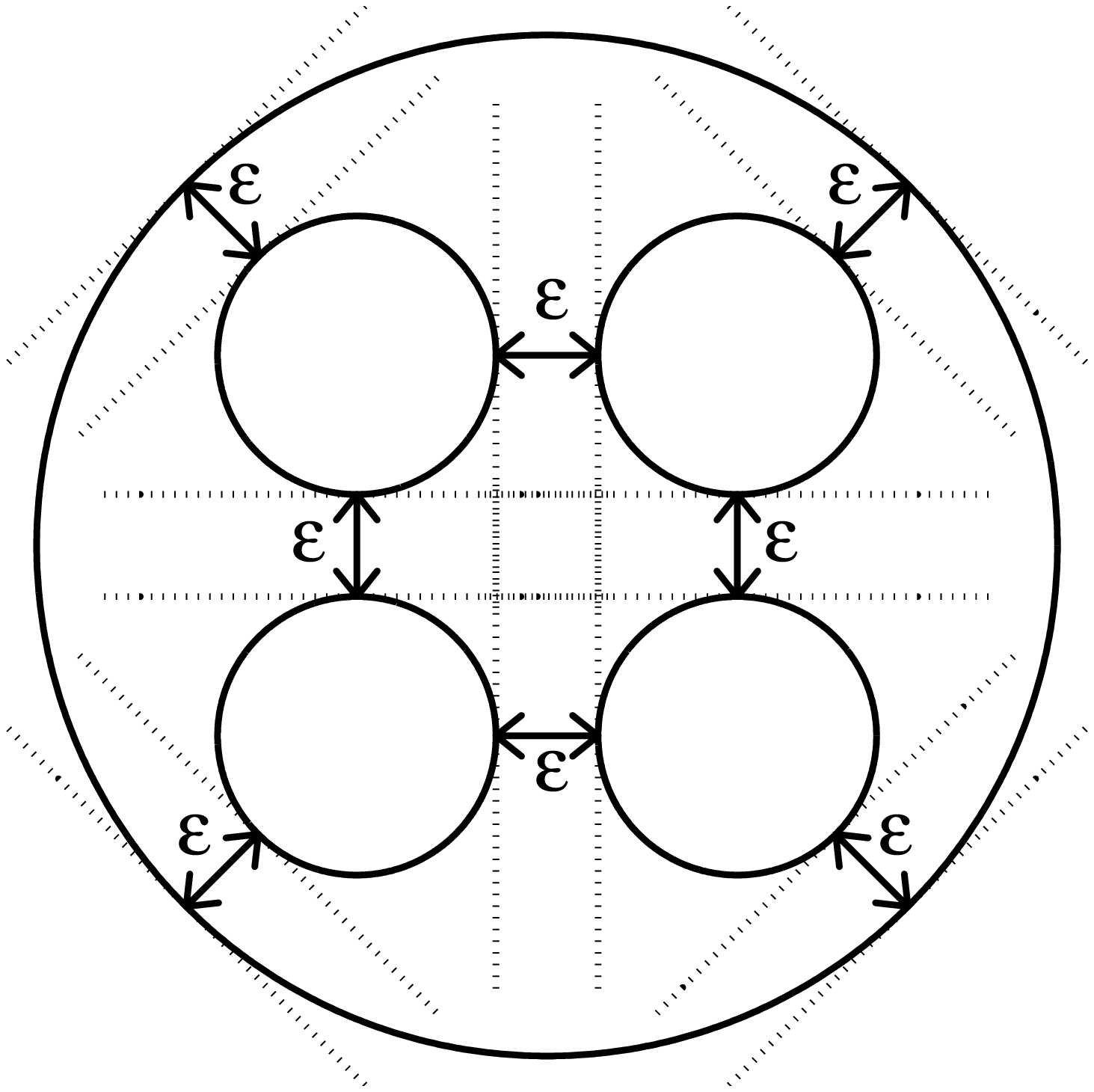}}}
\caption{\rm The domain of Examples~2 and~3.} 
\label{f:ex2-fig}
\end{figure}

\vspace{0.25cm}
{\noindent\bf Example 2.}
We consider a bounded multiply connected domains $G$ of connectivity $5$ ($m=4$). The boundary of $G$ consists of $5$ circles with variable distance $\varepsilon$ between these circles (see Fig.~\ref{f:ex2-fig}). The external circle is the unit circle. The internal four circles have the same radius 
\[
\frac{2-\varepsilon(2+2\sqrt{2})}{2+2\sqrt{2}}
\]
and the centres 
\[
\pm\left(\frac{2-\varepsilon}{2+2\sqrt{2}}\right)\pm\i\left(\frac{2-\varepsilon}{2+2\sqrt{2}}\right).
\]

In this example, we shall show the effect of the distance $\varepsilon$ on the accuracy of the functions \verb|FBIE| and \verb|FBIEad|. We shall consider the same functions $\gamma$, $\mu$, $h$, and $\phi$ as in~Example~1. For the piecewise constant function function $\theta$ in~(\ref{e:thet}), we shall consider two cases of the function $\theta(t)$. In the first case, we consider the constant function $\theta$
\[
\theta(t)=\pi/2\quad {\rm for\; all\;\; }t\in J,
\]
 i.e., $\theta$ has the same value on all boundary components. In the second case, we consider the non-constant function 
\[
\theta(t)=(\pi/2,0,\pi/2,0,\pi/2).
\]

\begin{figure}[!ht]%
\centerline{\scalebox{0.7}{\includegraphics{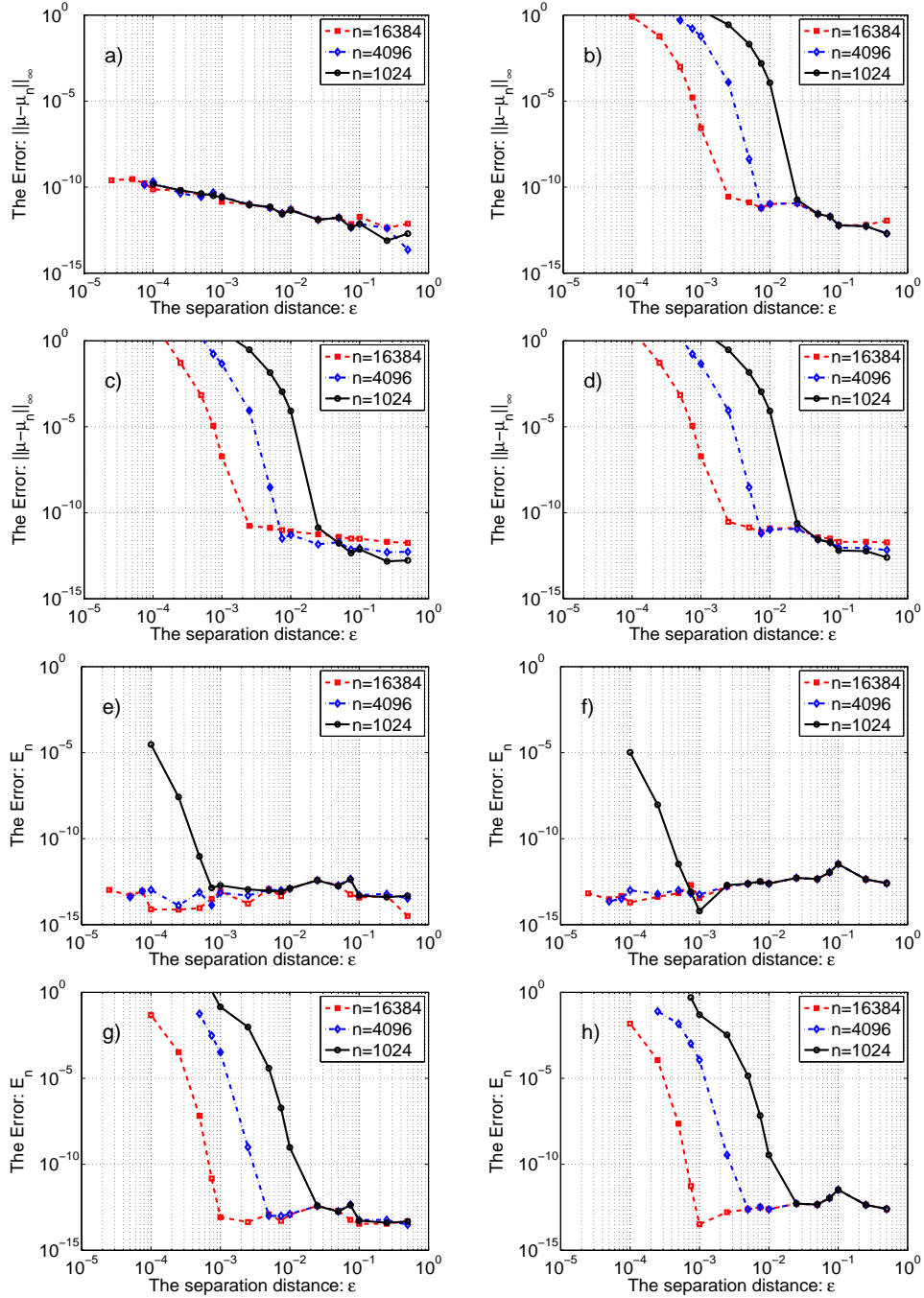}}}
\caption{\rm The maximum error norm $\|\mu-\mu_n\|_\infty$ and the error $E_n$ as a function of the separation distance $\varepsilon$. The figures on the left column obtained with the constant function $\theta$ and the figures on the right column obtained with the non-constant function $\theta$ using: (a,b) {\tt FBIE};  (c,d) the integral equation~(\ref{e:ie}) without singularity subtraction; (e,d) {\tt FBIEad}; (g,h) the integral equation~(\ref{e:ie*}) without singularity subtraction.} 
\label{f:ex2-err}
\end{figure}

\begin{figure}[!ht]%
\centerline{\scalebox{0.7}{\includegraphics{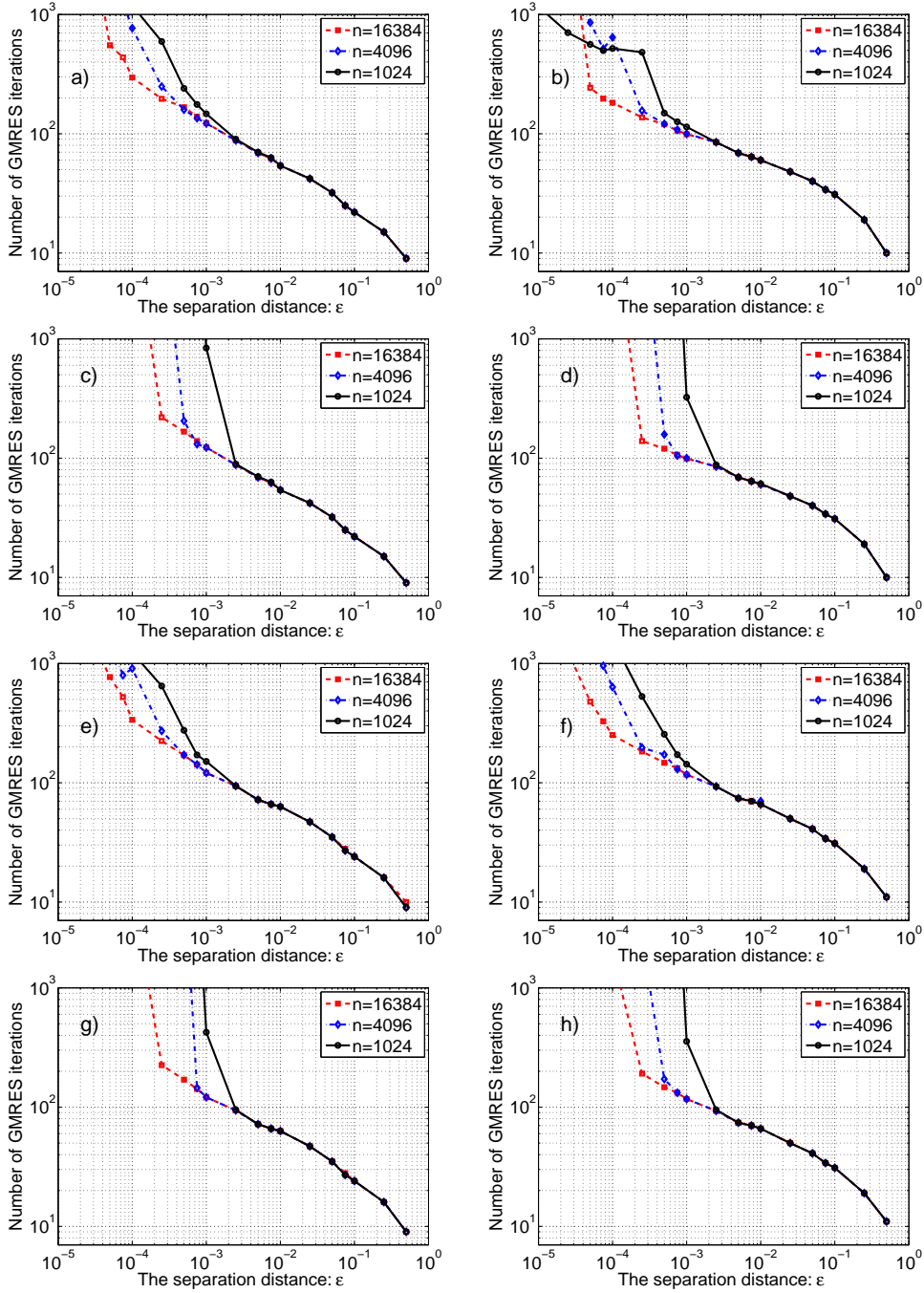}}}
\caption{\rm The number of GMRES iterations as a function of the separation distance $\varepsilon$. The figures on the left column obtained with the constant function $\theta$ and the figures on the right column obtained with the non-constant function $\theta$ for: (a,b) {\tt FBIE};  (c,d) the integral equation~(\ref{e:ie}) without singularity subtraction; (e,d) {\tt FBIEad}; (g,h) the integral equation~(\ref{e:ie*}) without singularity subtraction.} 
\label{f:ex2-iter}
\end{figure}

\begin{figure}[!ht]%
\centerline{\scalebox{0.7}{\includegraphics{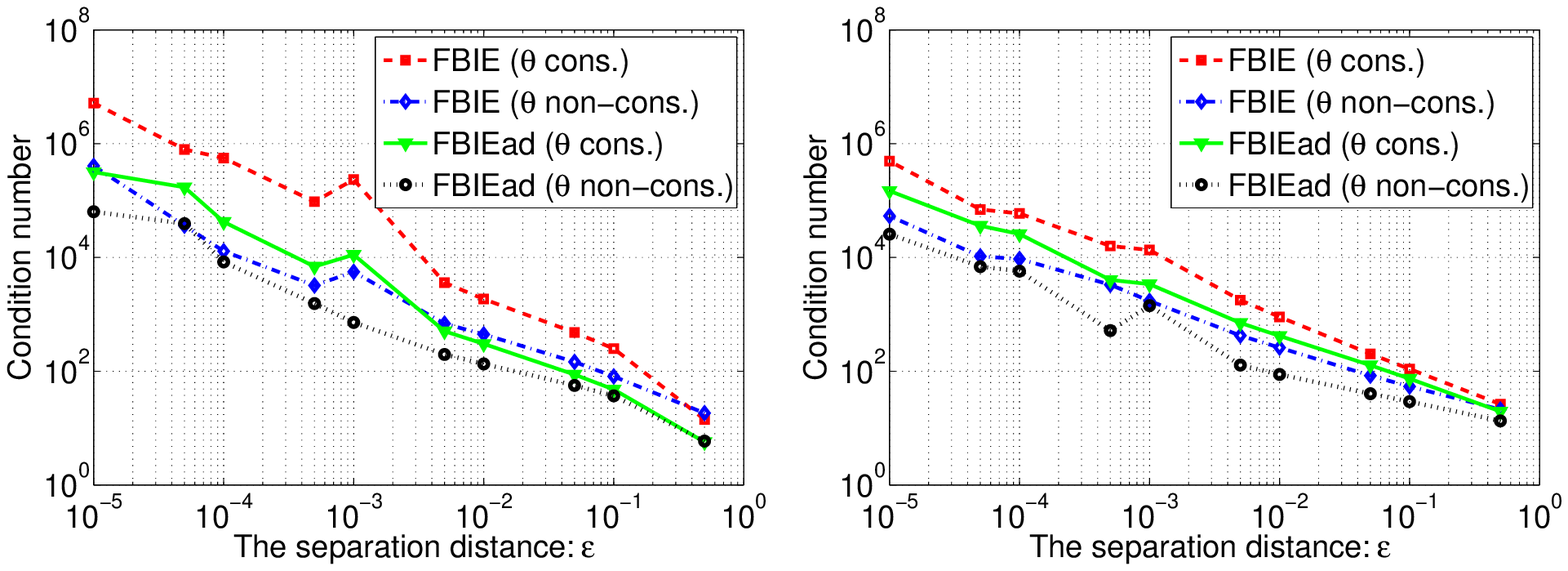}}}
\caption{\rm The condition number of the coefficient matrices of the linear systems for the bounded domain in Figure~\ref{f:ex2-fig} (left) and the unbounded domain in Figure~\ref{f:ex4-fig} (right).} 
\label{f:cond}
\end{figure}

The numerical results are shown in Figure~\ref{f:ex2-err}(a,b) for \verb|FBIE| and in Figure~\ref{f:ex2-err}(e,f) for \verb|FBIEad| where the error norms $\|\mu-\mu_n\|_\infty$ and $\|h-h_n\|_\infty$, and the error $E_n$ are defined and computed as in Example~1 with ${\tt iprec}=4$, ${\tt restart}=25$, ${\tt gmrestol}=10^{-12}$, and ${\tt maxit}=40$. To show the importance of the singularity subtraction in~(\ref{e:ie-m}) and~(\ref{e:ie*-m}), we solve the integral equations~(\ref{e:ie}) and~(\ref{e:ie*}) using the same method used in the functions \verb|FBIE| and \verb|FBIEad| but without singularity subtraction. The errors are shown in Figure~\ref{f:ex2-err}(c,d) for the integral equation~(\ref{e:ie}) and in Figure~\ref{f:ex2-err}(g,h) for the integral equation~(\ref{e:ie*}).
In general, the numerical results obtained with the functions \verb|FBIE| and \verb|FBIEad| are much better than the results obtained without singularity subtraction. If $\theta$ is constant, then the results obtained with \verb|FBIE| is much better than the results obtained without singularity subtraction (see Figure~\ref{f:ex2-err}(a,c)). The function \verb|FBIE| gives accurate results even for very small $\varepsilon$ when $n$ is sufficiently large. If $\theta$ is a non-constant function, then the accuracy of the results obtained by the function \verb|FBIE| is almost the same accuracy of the results obtained without singularity subtraction (see Figure~\ref{f:ex2-err}(b,d)). The function \verb|FBIEad| gives accurate results even for very small $\varepsilon$ when $n$ is sufficiently large for both the constant and the non-constant function $\theta$. For both cases, the results obtained with \verb|FBIEad| is much better than the results obtained without singularity subtraction (see Figure~\ref{f:ex2-err}(e,f,g,h)). 
For the functions \verb|FBIE| and \verb|FBIEad| as well as the methods without singularity subtraction, the number of GMRES iterations increase as $\varepsilon$ decreases (see Figure~\ref{f:ex2-iter}(a--h)). The condition number of the coefficient matrices of the linear systems for the functions \verb|FBIE| and \verb|FBIEad| is shown in Figure~\ref{f:cond} (left). The condition number is computed with $n=1024$ using the MATLAB condition number estimation function \verb|condest|. The condition number increase as $\varepsilon$ decreases. 

\vspace{0.25cm}
{\noindent\bf Example 3.}
We use the methods presented in~\cite{Nas-siam09,Nas-jmaa11} which is based on the integral equation with the generalized Neumann kernel~(\ref{e:ie}) and the method presented in~\cite{Nas-inv,Yun-sp,Yun-inv} which is based on the integral equation with the adjoint generalized Neumann kernel~(\ref{e:ie*}) to compute the conformal mapping from the bounded multiply connected domains $G$ of Example~2 (see Fig.~\ref{f:ex2-fig}) onto the disc with circular slits and the disc with both circular and radial slits. The integral equations are solved using the functions \verb|FBIE| and \verb|FBIEad|  with ${\tt iprec}=4$, ${\tt restart}=25$, ${\tt gmrestol}=10^{-12}$, and ${\tt maxit}=40$. For the disc with circular slits, the function $\theta$ is a constant function where $\theta(t)=\pi/2$ for all $t\in J$. For the disc with both circular and radial slits, the function $\theta$ is not a constant function. Its values is $\pi/2$ on the external boundary, $\pi/2$ on the boundaries which mapped to circular slits, and $0$ on the boundaries which mapped to radial slits. In this example, for the disc with both circular and radial slits, we shall assume that $\theta(t)=(\pi/2,0,\pi/2,0,\pi/2)$.

\begin{figure}[!ht]%
\centerline{\scalebox{0.8}{\includegraphics{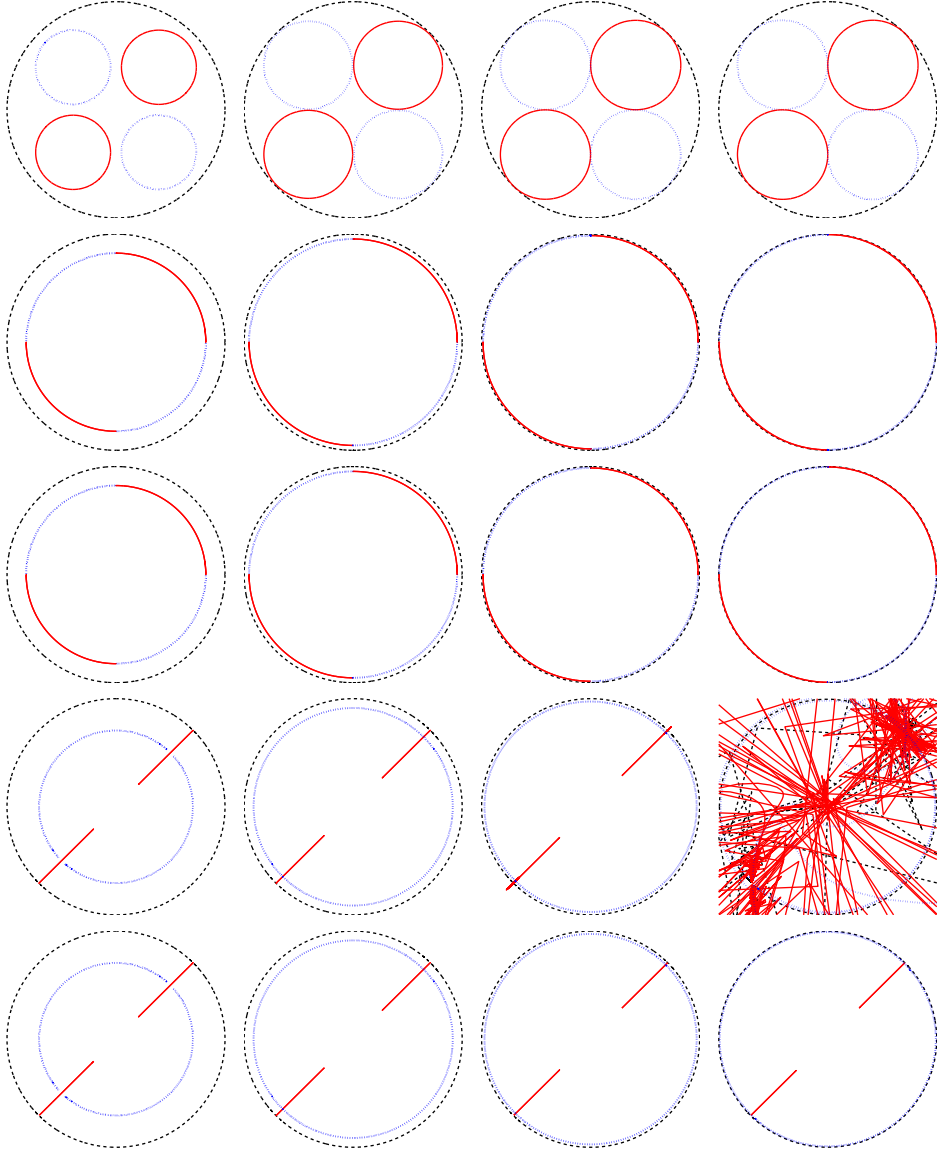}}}
\caption{\rm In the first row, the domain $G$ for separation distance $\varepsilon=10^{-1}$, $10^{-2}$, $10^{-3}$, $10^{-4}$. In the second row, the image of $G$ obtained with {\tt FBIE} for the constant function $\theta$. In the third row, the image of $G$ obtained with {\tt FBIEad} for the constant function $\theta$. In the fourth row, the image of $G$ obtained with {\tt FBIE} for the non-constant function $\theta$. In the fifth row, the image of $G$ obtained with {\tt FBIEad} for the non-constant function $\theta$.} 
\label{f:ex3}
\end{figure}

\begin{figure}[!ht]%
\centerline{\scalebox{0.7}{\includegraphics{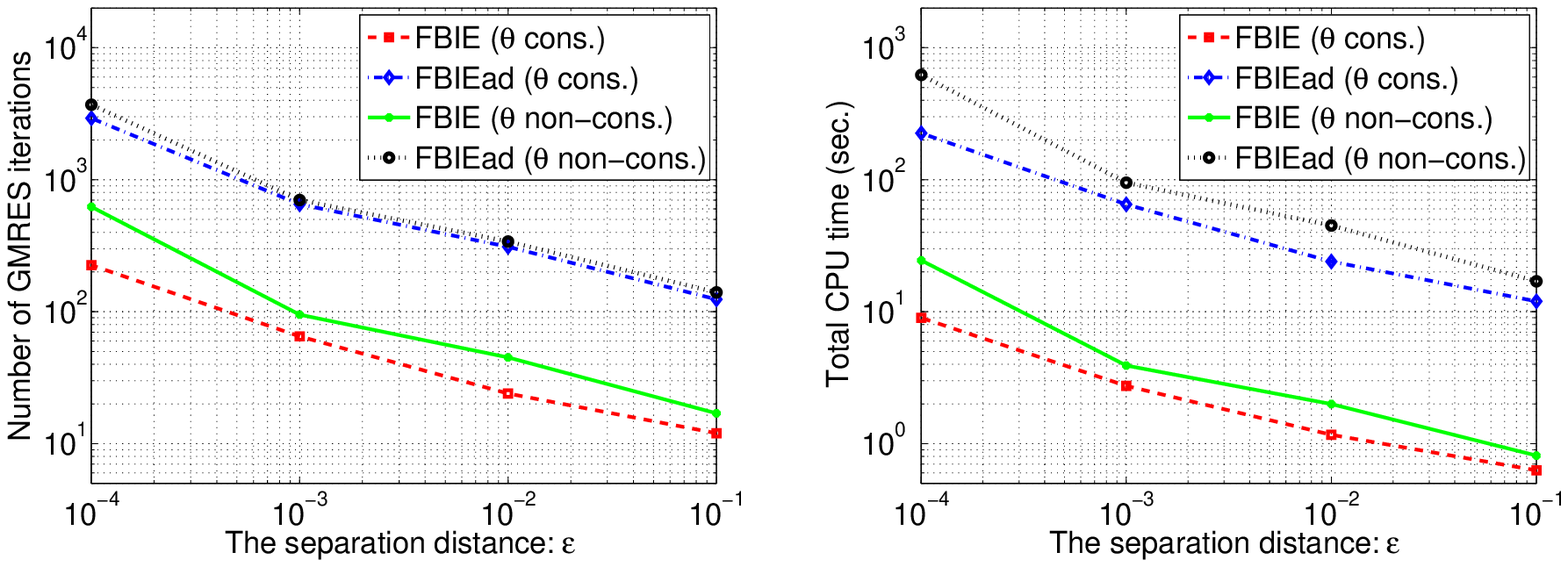}}}
\caption{\rm The number of GMRES iterations and the total CPU time vs. the separation distance $\varepsilon$.} 
\label{f:ex3-err}
\end{figure}

The original domain $G$ for separation distance $\varepsilon=10^{-1},\,10^{-2},\,10^{-3},\,10^{-4}$ is shown in Figure~\ref{f:ex3}. The methods presented in~\cite{Nas-siam09,Nas-jmaa11} is based on solving the integral equation~(\ref{e:ie}) then computing the function $h$ by~(\ref{e:h}) to obtain the boundary values of the mapping function. Thus the complexity of the method based on the integral equation with the generalized Neumann kernel~(\ref{e:ie}) is $O(m+1)n\ln n)$ (see~\cite{Nas-siam13}). The images of the original domains obtained with \verb|FBIE| for $n=4096$ are shown in the second row of Figure~\ref{f:ex3} for constant $\theta$ and in the fourth row of Figure~\ref{f:ex3} for non-constant $\theta$. The number of GMRES iterations and the total CPU time (in seconds) required to obtain the boundary values of the mapping function for $n=4096$ vs. the separation distance $\varepsilon$ are shown in Figure~\ref{f:ex3-err}.   

For the method presented in~\cite{Nas-inv,Yun-sp,Yun-inv}, it is based solving $m+2$ integral equation with the adjoint generalized Neumann kernel of the form~(\ref{e:ie*}) to compute the boundary values of the mapping function. We need to solve $m+1$ integral equations to obtain the parameters of the canonical domain and one integral equation to obtain the derivative of the boundary correspondence function. The complexity of solving these $m+2$ integral equations using the function \verb|FBIEad| is $O((m+2)(m+1)n)$. By obtaining the derivative of the boundary correspondence function, we need to use the FFT in each boundary component $J_j$, $j=0,1,\ldots,m$, to obtain the boundary correspondence function. The complexity of computing the boundary correspondence function from its derivative is $O((m+1)n\ln n)$. Thus, the complexity of the method based on the integral equation with the adjoint generalized Neumann kernel~(\ref{e:ie*}) is $O((m+1)(m+2+\ln n)n)$. The images of the original domains obtained with \verb|FBIEad| for $n=4096$ are shown in the third row of Figure~\ref{f:ex3} for constant $\theta$ and in the fifth row of Figure~\ref{f:ex3} for non-constant $\theta$. The total number of GMRES iterations (for solving the $6$ integral equations) and the total CPU time (in seconds) required to obtain the boundary values of the mapping function for $n=4096$ vs. the separation distance $\varepsilon$ is shown in Figure~\ref{f:ex3-err}.   

As explained in Example~2, the function \verb|FBIE| gives accurate results for the constant function $\theta$. For the non-constant function $\theta$, we get accurate results for $\varepsilon=10^{-1},10^{-2}$. For $\varepsilon=10^{-3}$ the radial slits goes outside of the unit disc and for $\varepsilon=10^{-4}$ the obtained figure is incorrect (see the third and fourth figures in the fourth row in Figure~\ref{f:ex3}). The function \verb|FBIEad| gives accurate results for both the constant function and non-constant function $\theta$. However, the complexity of the method based on the integral equation with the adjoint generalized Neumann kernel~(\ref{e:ie*}) is larger than the complexity of the method based on the integral equation with the generalized Neumann kernel~(\ref{e:ie}) specially for large~$m$. As Figure~\ref{f:ex3-err} shows, the number of GMRES iterations as well as the CPU time of the method based on the integral equation with the adjoint generalized Neumann kernel~(\ref{e:ie*}) is much larger than the the number of GMRES iterations and the CPU time of the method based on the integral equation with the generalized Neumann kernel~(\ref{e:ie}).

\begin{figure}[!ht]%
\centerline{\scalebox{0.3}{\includegraphics{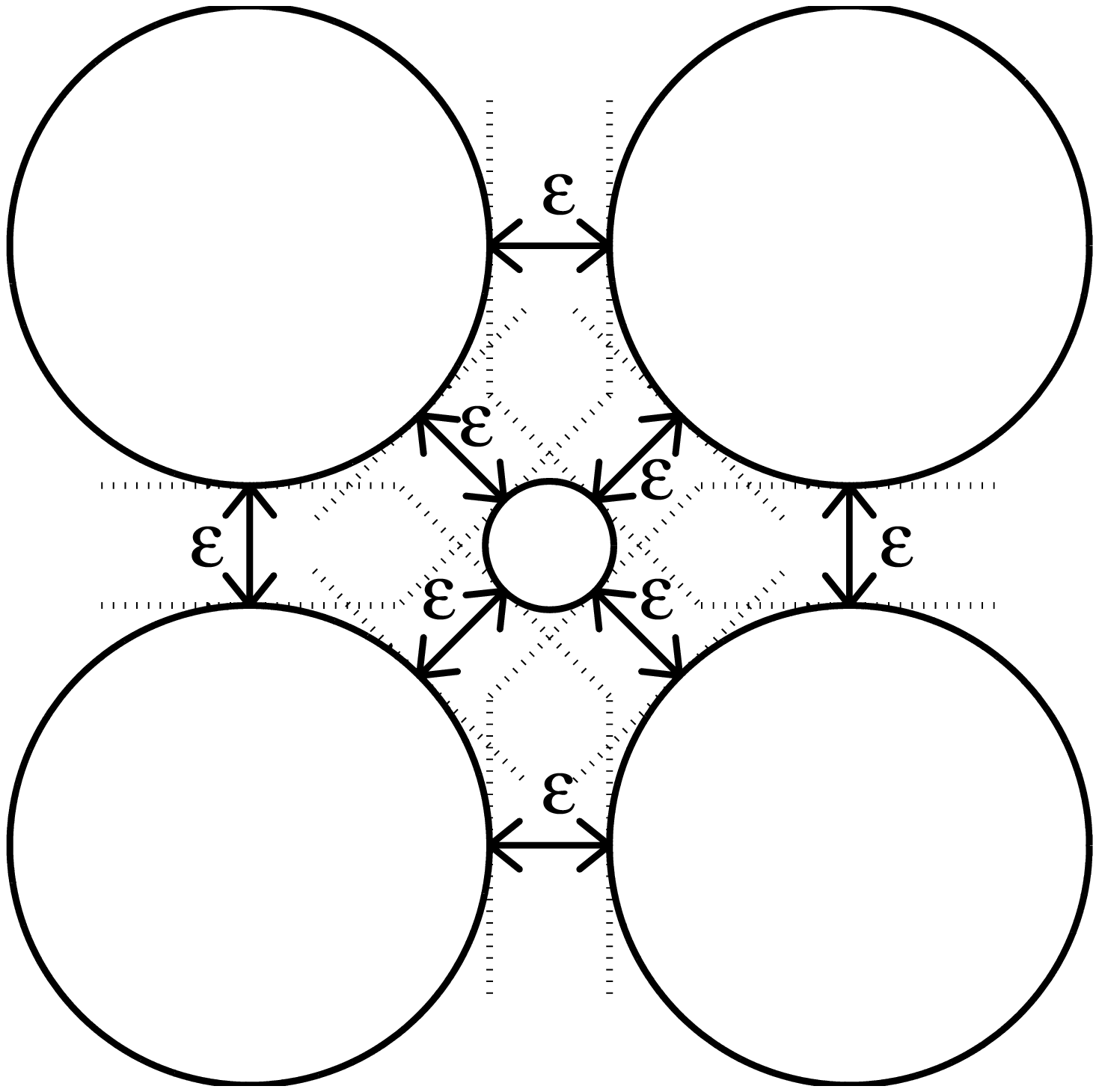}}}
\caption{\rm The domain of Example~4.} 
\label{f:ex4-fig}
\end{figure}

\begin{figure}[!ht]%
\centerline{\scalebox{0.8}{\includegraphics{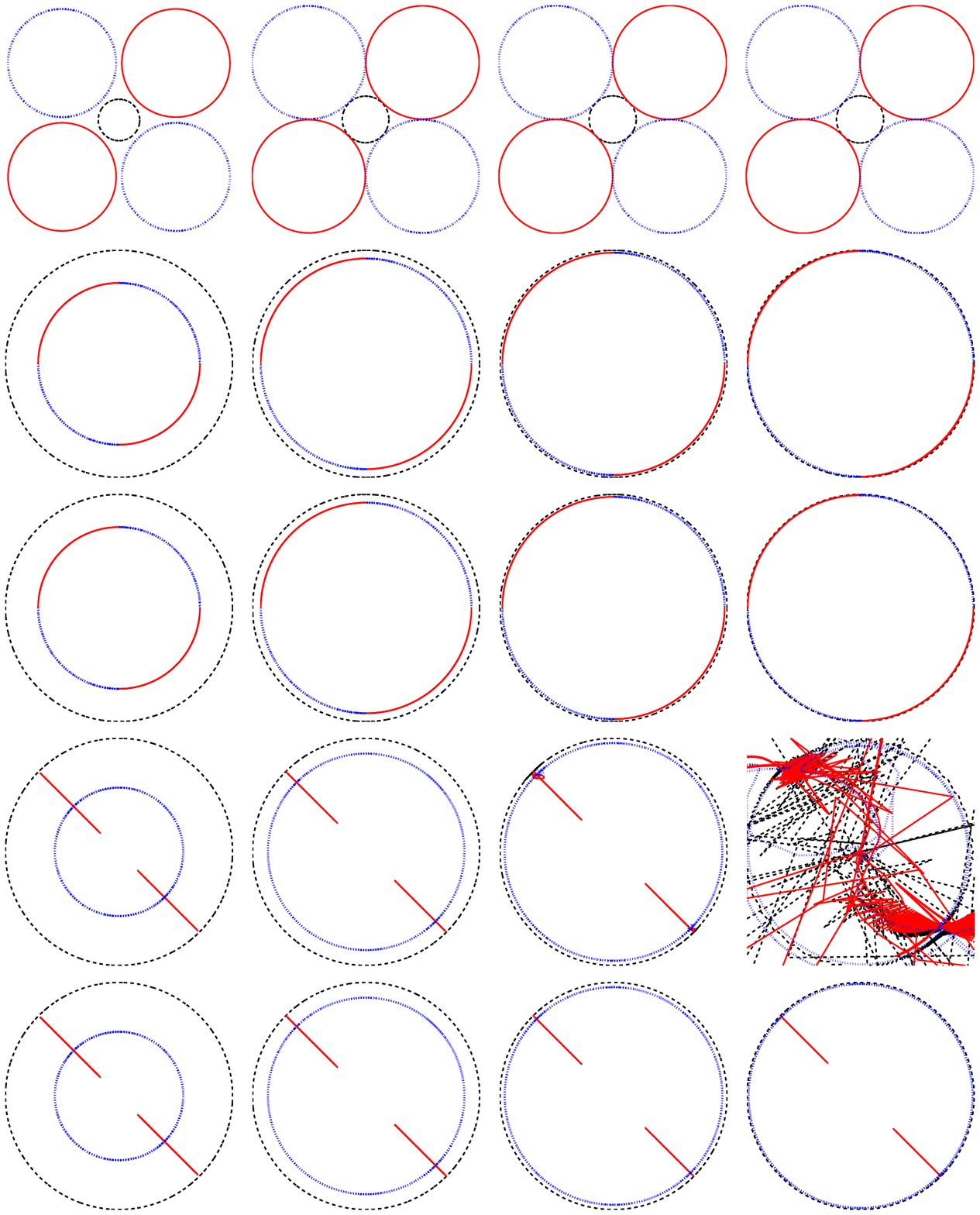}}}
\caption{\rm In the first row, the domain $G$ for separation distance $\varepsilon=10^{-1}$, $10^{-2}$, $10^{-3}$, $10^{-4}$. In the second row, the image of $G$ obtained with {\tt FBIE} for the constant function $\theta$. In the third row, the image of $G$ obtained with {\tt FBIEad} for the constant function $\theta$. In the fourth row, the image of $G$ obtained with {\tt FBIE} for the non-constant function $\theta$. In the fifth row, the image of $G$ obtained with {\tt FBIEad} for the non-constant function $\theta$.} 
\label{f:ex4}
\end{figure}

\begin{figure}[ht]
\centerline{\scalebox{0.7}{\includegraphics{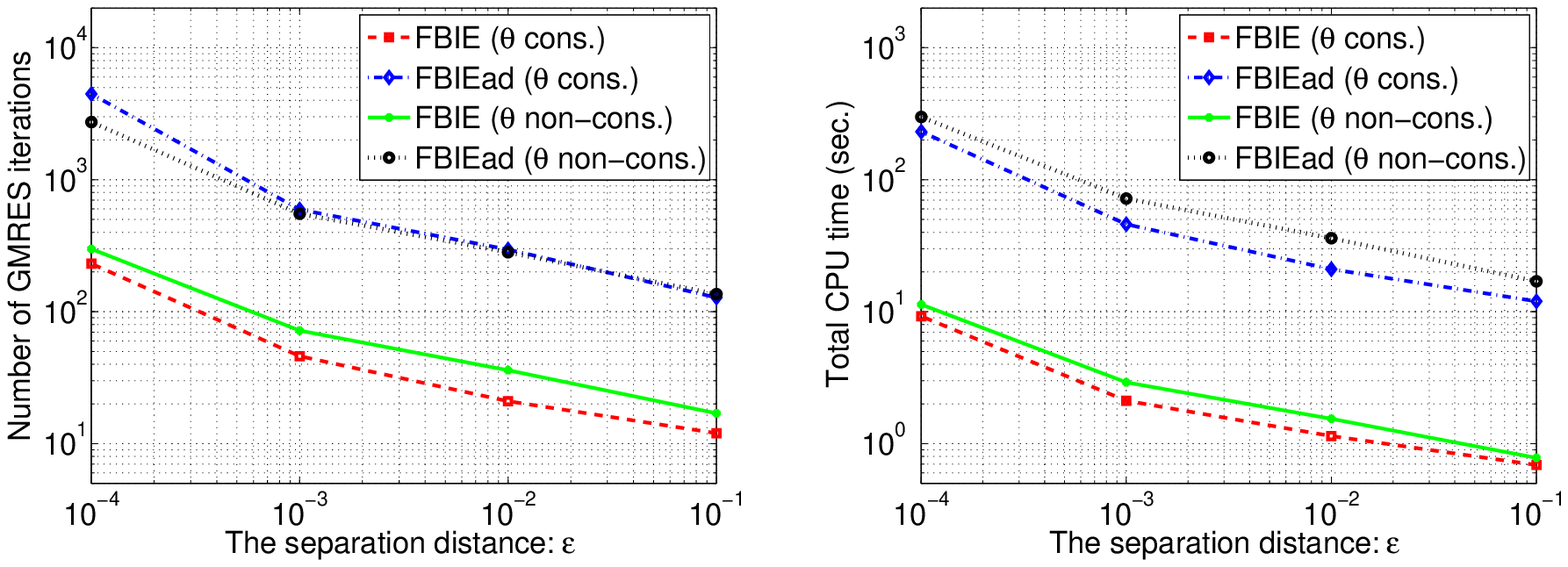}}}
\caption{\rm The number of GMRES iterations and the total CPU time vs. the separation distance $\varepsilon$.} 
\label{f:ex4-err}
\end{figure}

\vspace{0.25cm}
{\noindent\bf Example 4.}
We repeat the Example~3 for an unbounded multiply connected domains $G$ of connectivity $5$ ($m=4$). The boundary of $G$ consists of $5$ circles where the distance between any two circles is $\varepsilon$ (see Fig.~\ref{f:ex4-fig}). The centre of the circle in the centre is $0$ and its radius is $\sqrt{2}-1-\frac{\varepsilon}{2}$. The other four circles have the same radius $1-\frac{\varepsilon}{2}$ and the centres $\pm1\pm\i$. The original domain $G$ for separation distance $\varepsilon=10^{-1},\,10^{-2},\,10^{-3},\,10^{-4}$ and its images computed with $n=4096$ are shown in Figure~\ref{f:ex4}. The total number of GMRES iterations and the total CPU time (in seconds) required to obtain the boundary values of the mapping function for $n=4096$ vs. the separation distance $\varepsilon$ are shown in Figure~\ref{f:ex4-err}. The condition number of the coefficient matrices of the linear systems for the functions \verb|FBIE| and \verb|FBIEad| computed with $n=1024$ using the MATLAB condition number estimation function \verb|condest| is shown in Figure~\ref{f:cond} (right).

\section{Conclusions}

This paper presented two new numerical methods for fast computing of the solutions of the boundary integral equations with the generalized Neumann kernel and the adjoint generalized Neumann kernel. The methods are based on discretizing the integral equations using the Nystr\"om method with the trapezoidal rule then solving the obtained linear systems using the GMRES method combined with the FMM. The complexity of the presented methods is $O((m+1)n\ln n)$ for the integral equations with the generalized Neumann kernel and $O((m+1)n)$ for the integral equations with the adjoint generalized Neumann kernel. The presented methods are fast, accurate, and can be used for domains with high connectivity, complex geometry, and close boundaries. 

Based on the presented methods, two MATLAB functions \verb|FBIE| and \verb|FBIEad| are presented for fast computing of the solutions of the integral equations with the generalized Neumann kernel and the adjoint generalized Neumann kernel, respectively. The accuracy of the presented functions \verb|FBIE| and \verb|FBIEad| for domains with close boundaries has been studied in Examples~2, 3, and~4. For both functions, the number of GMRES iterations, the total CPU time, and the condition number of the coefficient matrices increase as the distance $\varepsilon$ decreases. The singularity subtraction increases the accuracy significantly. For the function \verb|FBIEad|, the accuracy of the function \verb|FBIEad|, the number of GMRES iterations, the total CPU time, and the condition number of the coefficient matrices are not affected by $\theta$ being constant function or not. 
For the function \verb|FBIE|, on one hand, the accuracy of the function \verb|FBIE| is affected by $\theta$ being constant function or not. We get accurate results for constant function $\theta$ even for very small distance $\varepsilon$. For the non-constant function $\theta$, we get accurate results for moderate small $\varepsilon$ then the accuracy is getting worse as the distance $\varepsilon$ becomes very small. On the other hands, the number of GMRES iterations, the total CPU time and the condition number of the coefficient matrices are not affected by $\theta$ being constant function or not. 
A possible reason for the effect of $\theta$ being constant function or not on the accuracy of \verb|FBIE| is that the function \verb|FBIE| is based on discretizing a Fredholm integral operator $\bN$ and a singular integral operator $\bM$. The function \verb|FBIEad| is based on discretizing only a Fredholm integral operator $\bN^\ast$. 
The number of GMRES iterations, the total CPU time, and the condition number of the coefficient matrices, which are not affected by $\theta$, depends only the Fredholm integral operators $\bN$ and $\bN^\ast$. So, the discretization of the operator $\bM$ could be the reason behind the effect of $\theta$ on the accuracy of \verb|FBIE|. There is numerical evidence that the accuracy of the function \verb|FBIE| is affected by the function $\theta$ if the operator $\bM$ is discretized by the method presented in~\cite{Nas-siam13}.


\end{document}